\documentclass[final,leqno,onefignum,onetabnum]{siamltex1213}

\usepackage{graphicx} % for pdf, bitmapped graphics files
\usepackage{mathptmx} % assumes new font selection scheme installed
\usepackage{times} % assumes new font selection scheme installed
\usepackage{amsmath} % assumes amsmath package installed
\usepackage{amssymb}  % assumes amsmath package installed
\usepackage{color}
\usepackage{epstopdf}
\usepackage{subfigure}
\usepackage{tikz}
\usetikzlibrary{matrix,arrows,shapes.geometric}
\usepackage[all,cmtip]{xy}

\def\a {\alpha}
\def\b {\beta}
\def\e {\varepsilon}
\def\s {\sigma}
\def\w {\omega}

\def\g {\gamma}
\def\l {\lambda}

\newtheorem{remark}{Remark}
\newtheorem{example}{Example}

\title{Fixed-Endpoint Optimal Control of Bilinear Ensemble Systems\thanks{This work was
supported in part by the National Science Foundation under the awards CMMI-1301148, CMMI-1462796, and ECCS-1509342.}} 

\author{Shuo Wang\thanks{Department of Electrical and Systems Engineering, Washington University, St. Louis, Missouri, 63130, USA
(\email{swang35@wustl.edu}). } \and Jr-Shin Li\thanks{Department of Electrical and Systems Engineering, Washington University, St. Louis, Missouri, 63130, USA
(\email{jsli@wustl.edu}). Questions, comments, or corrections
to this document may be directed to this email address.}}

\begin{document}
\maketitle
\slugger{mms}{xxxx}{xx}{x}{x--x}%slugger should be set to mms, siap, sicomp, sicon, sidma, sima, simax, sinum, siopt, sisc, or sirev

\begin{abstract}
Optimal control of bilinear systems has been a well-studied subject in the area of mathematical control. However, techniques for solving emerging optimal control problems involving an ensemble of structurally identical bilinear systems are underdeveloped. In this work, we develop an iterative method to effectively and systematically solve these challenging optimal ensemble control problems, in which the bilinear ensemble system is represented as a time-varying linear ensemble system at each iteration and the optimal ensemble control law is then obtained by the singular value expansion of the input-to-state operator that describes the dynamics of the linear ensemble system. We examine the convergence of the developed iterative procedure and pose optimality conditions for the convergent solution. We also provide examples of practical control designs in magnetic resonance to demonstrate the applicability and robustness of the developed iterative method.
% obtained using a singular-value-decomposition based computational algorithm designed for solving optimal control problems involving linear ensemble systems. We examine the convergence of the developed iterative procedure and construct optimality conditions for the convergent solution. We also provide examples of practical control designs in magnetic resonance to demonstrate the applicability and robustness of the developed iterative method.
\end{abstract}


\begin{AMS}\end{AMS}

\pagestyle{myheadings}
\thispagestyle{plain}
\markboth{S. WANG AND J.-S. LI}{FIXED-ENDPOINT CONTROL FOR BILINEAR ENSEMBLES}

%%%%%%%%%%%%%%%%%%%%%%%%%%%%%%%
\section{Introduction}
Newly emerging fields in science and engineering, such as systems neuroscience, synchronization engineering, and quantum science and technology, give rise to new classes of optimal control problems that involve underactuated manipulation of individual and collective behavior of dynamic units in a large ensemble. Representative examples include neural stimulation for alleviating the symptoms of neurological disorders such as Parkinson's disease, where a population of neurons in the brain is affected by a small number of electrodes \cite{Ching2013a}; %one or two electrodes \cite{Ching2013a}; 
pulse designs for exciting and transporting quantum systems between desired states, where an ensemble of quantum systems is driven by a single or multiple pulses in a pulse sequence \cite{Brent06, Li_PNAS11}; and the engineering of dynamical structures for complex oscillator networks, where 
sequential patterns of a network of nonlinear rhythmic elements are created and altered by a mild global waveform \cite{Kiss:2007:1886-1889}.
%is designed to alter sequential patterns of nonlinear rhythmic elements \cite{Kiss:2007:1886-1889}. 
Solving these nontraditional and large-scale underactuated control problems requires the development of systematic and computationally tractable and effective methods.

Among these emerging control problems, in this paper, we will study fixed-endpoint optimal control problems involving bilinear ensemble systems, which arise from the domain of quantum control \cite{Chen14} and appear in a variety of other different fields, such as cancer chemotherapy \cite{Schaettler2002} and robotics \cite{Becker12}. The control of bilinear systems has been a well-studied subject in the area of mathematical control. From Pontryagin's maximum principle to spectral collocation methods, a wide variety of theoretical and computational methods have been developed to solve optimal control problems of bilinear systems \cite{Aganovic1994, Mohler2000}. In particular, the numerical methods are in principle categorized into direct, e.g., pseudospectral methods \cite{Gong06,Ross04}, and indirect approaches, e.g., indirect transcription method \cite{Rao09} and %indirect 
shooting methods \cite{Bertolazzi2005}. Implementing these existing numerical methods to solve optimal control problems involving an ensemble, i.e., a large number (finitely or infinitely many) or a parameterized family, of bilinear systems %with slightly different dynamics characterized by common system parameters 
may encounter low efficiency, slow convergence, and instability issues, because most of these methods rely on suitable discretization of the continuous-time dynamics into a large-scale nonlinear program (LSNLP). In addition, the global constraint for such an optimal ensemble control problem, in which each individual system receives the same control input, makes the discretized LSNLP very restrictive and intractable to solve or even to find a feasible solution \cite{Li_TAC12_QCP}.

On the other hand, optimal control problems involving a linear system, or a linear ensemble system, are often computationally tractable and analytically solvable for many special cases, such as the linear quadratic regulator (LQR) \cite{Brockett70} and the minimum-energy control of harmonic oscillator ensembles \cite{Li_TAC11}. This suggests a bypass to solve optimal control problems of bilinear ensemble systems through solving that of linear ensemble systems and motivates the development of the iterative method in this work. %, for which an iterative method is developed in this work. 
The central idea is to represent the bilinear ensemble system as a linear ensemble system at each iteration, and then feasibly calculate the optimal control and trajectory for each iteration until a convergent solution is found. Iterative methods have been introduced and adopted to deal with diverse control design problems, including the free-endpoint quadratic optimal control of bilinear systems \cite{Hofer1988} and optimal state tracking for nonlinear systems \cite{Cimen2004}, while the fixed-endpoint problems along with the emerging problems that involve controlling a bilinear ensemble system remain unexplored. 

In this paper, we combine the idea of the aforementioned iterative method with our previous work on optimal control of linear ensemble systems to construct an iterative algorithm for solving optimal control problems involving a time-invariant bilinear ensemble system of the form,
\begin{eqnarray*}
	%\label{eq:bilinear_ensemble}
	\frac{d}{dt}{X(t,\b)}=A(\b)X(t,\b)+B(\b)u(t)+\Big(\sum_{i=1}^m u_i(t) B_i(\b)\Big) X(t,\b),
\end{eqnarray*}
where $X=(x_1,\ldots,x_n)^T\in M\subset\mathbb{R}^n$ denotes the state, $\b\in K\subset\mathbb{R}^d$ with $K$ compact and $d$ a positive integer, $u(t)=(u_1(t),\ldots,u_m(t))^T \in\mathbb{R}^m$ is the control, and the matrices $A(\b)\in\mathbb{R}^{n\times n}$, $B(\b)\in\mathbb{R}^{n\times m}$, and $B_i(\b)\in\mathbb{R}^{n\times n}$, $i=1,\ldots,m$, for $\b\in K$.

This paper is structured as follows. In the next section, we present the developed iterative method for fixed-endpoint optimal control of a time-invariant bilinear system, where we introduce a sweep method that accounts for %that guarantees the fulfillment of 
the terminal condition based on the notion of flow mapping from the optimal control theory. % through a set of different equations. 
In Section \ref{sec:convergence}, we examine the convergence of the iterative method using the fixed-point theorem. In Section \ref{sec:optimality}, we propose the conditions for global optimality of the convergent solution. Then, in Section \ref{sec:ensemble}, we extend the developed iterative method to solve optimal control problems involving bilinear ensemble systems and show the convergence of the method. Finally, examples and simulations of practical control design problems are illustrated in Section \ref{sec:examples} to demonstrate the applicability and robustness of the developed iterative procedure.

%%%%%%%%%%%%%%%%%%%%%%%%%%%%%%%%%%%%%%%%%%%%%%%%%%%%%%%%%%%%
\section{Iterative method for optimal control of bilinear systems}
\label{sec:single_bilinear}
We start with considering a fixed-endpoint, finite-time, quadratic optimal control problem involving a time-invariant bilinear system of the form
\begin{align}
	\min &\quad J=\frac{1}{2}\int_0^{t_f} \Big[x^T(t)Qx(t)+u^T(t)Ru(t)\Big] \, dt, \nonumber\\
	\label{eq:oc1}
	{\rm s.t.} &\quad \dot{x}=Ax+Bu+\Big[\sum_{i=1}^m u_i B_i\Big]x, \tag{P1}\\ %\Big(\sum_{i=1}^n x_j N_j\Big) u, \\
	&\quad x(0)=x_0, \quad x(t_f)=x_f, \nonumber
\end{align}
where $x(t)\in\mathbb{R}^n$ is the state and $u(t)\in\mathbb{R}^m$ is the control; $A\in\mathbb{R}^{n\times n}$, $B_i\in\mathbb{R}^{n\times n}$, and $B\in\mathbb{R}^{n\times m}$ are constant matrices; $R\in\mathbb{R}^{m\times m}\succ 0$ is positive definite and $Q\in\mathbb{R}^{n\times n}\succeq 0$ is positive semi-definite; and $x_0,x_f\in\mathbb{R}^n$ are the initial and the desired terminal state, respectively. We first represent the time-invariant bilinear system in \eqref{eq:oc1} as a time-varying linear system,
\begin{eqnarray}
	\label{eq:linear}
	\dot{x}(t)=Ax+Bu+\left[\sum_{j=1}^n x_j(t) N_j\right]u,
\end{eqnarray}
in which we write the bilinear term $\big(\sum_{i=1}^m u_iB_i\big)x=\big(\sum_{j=1}^n x_jN_j\big) u$ with $x_j$ the $j^{th}$ element of $x$, $N_j\in\mathbb{R}^{n\times m}$ for $j=i,\dots,n$, and $u=(u_1,\ldots,u_m)^T\in\mathbb{R}^m$. Then, we solve this optimal control problem by Pontryagin's maximum principle. The Hamiltonian of this problem is
\begin{equation}
H(x,u,\lambda)=\frac{1}{2}(x^TQx+u^TRu)+\lambda^T \Big\{Ax+\Big[B+(\sum_{j=1}^n x_jN_j)\Big] u\Big\},
\end{equation} 
where $\lambda(t)\in\mathbb{R}^n$ is the co-state vector. The optimal control is then obtained by the necessary condition, $\frac{\partial H}{\partial u}=0$, given by
\begin{equation}
	\label{eq:u}
	u^* = -R^{-1}\Big( B+\sum_{j=1}^n x_jN_j\Big)^T\lambda,
\end{equation}
and the optimal trajectory of the state $x$ and the co-state $\lambda$ satisfy, for $t\in[0,t_f]$,
\begin{align}
	\label{eq:state}
	\dot{x}_i & = \big[Ax\big]_i-\Big[\big(B+\sum_{j=1}^n x_jN_j\big)R^{-1}\big(B+\sum_{j=1}^n x_jN_j\big)^T\lambda\Big]_i, \\
	\label{eq:costate}
	\dot{\lambda}_i & = -\big[Qx\big]_i-\big[A^T\lambda\big]_i+\lambda^T\Big\{N_iR^{-1}\big(B+\sum_{j=1}^n x_jN_j\big)^T+\big(B+\sum_{j=1}^n x_jN_j\big)R^{-1}N_i^T\Big\}\lambda
\end{align}
with the boundary conditions $x(0)=x_0$ and $x(t_f)=x_f$, where $x_i$, $\l_i$ and $[\,\cdot\,]_i$, $i=1,\dots,n$, are the $i^{th}$ component of the associated vectors. By the following change of variables,
\begin{align}
	\label{eq:A}
	& \tilde{A}_{ij} = A_{ij}-\Big[(N_j R^{-1} \big(B+\sum_{j=1}^n x_jN_j\big)^T+\big(B+\sum_{j=1}^n x_jN_j) R^{-1} N_j^T\big)\lambda\Big]_i, \\
	\label{eq:B}
	&\tilde{B}R^{-1}\tilde{B}^T = BR^{-1}B^T-\big(\sum_{j=1}^n x_jN_j\big)R^{-1}\big(\sum_{j=1}^n x_jN_j\big)^T,\\
	\label{eq:Q}
	&\tilde{Q}=Q,
\end{align}
we can rewrite \eqref{eq:state} and \eqref{eq:costate} into the form
\begin{align}
	\label{eq:state1}
	\dot{x} &= \tilde{A}x-\tilde{B}R^{-1}\tilde{B}^T\lambda, \quad x(0)=x_0, \quad  x(t_f) = x_f, \\ 
	\label{eq:costate1}
	\dot{\lambda} &= -\tilde{Q}x-\tilde{A}^T\lambda,
\end{align}
which coincides with the canonical form of the state and co-state equations characterizing the optimal trajectories for the analogous optimal control problem involving the time-invariant linear system $\dot{x}=\tilde{A}x+\tilde{B}u$ \cite{Schaettler13}. In this way, the optimal state and co-state trajectories for the optimal control problem \eqref{eq:oc1} involving a time-invariant bilinear system %as presented in \eqref{eq:oc1} 
are now expressed in terms of the equations related to a time-varying linear system as in \eqref{eq:state1} and \eqref{eq:costate1}.

Using this ``linear-system representation'' together with the Sweep method \cite{Schaettler13, Bryson75}, we will solve the optimal control problem \eqref{eq:oc1} in an iterative manner. Specifically, we will consider at each iteration the fixed-endpoint linear quadratic optimal control problem,
\begin{align}
	\min &\quad J=\frac{1}{2}\int_0^{t_f} \Big[(x^{(k+1)})^T(t)Qx^{(k+1)}(t)+(u^{(k+1)})^T(t)Ru^{(k+1)}(t)\Big] \, dt, \nonumber\\
	\label{eq:oc2}
	{\rm s.t.} &\quad %\dot{x}^{(k+1)}(t)=A^{(k)}x^{(k+1)}+B^{(k)} u^{(k+1)}+\left[\sum_{j=1}^n x_j^{(k)}(t) N_j\right] u^{(k+1)}, \tag{P2}\\
	\dot{x}^{(k+1)}(t)=\tilde{A}^{(k)} x^{(k+1)}+\tilde{B}^{(k)} u^{(k+1)} \tag{P2}\\  
	&\quad x^{(k+1)}(0)=x_0, \quad x^{(k+1)}(t_f)=x_f, \nonumber
\end{align}
by treating the previous trajectory $x^{(k)}$ as a known quantity, where $k\in\mathbb{N}$ denotes the iteration. In the following sections, we will introduce the Sweep method and present the iterative procedure.

% =======================================================
\subsection{Sweep method for fixed-endpoint problems} 
\label{sec:sweep}
Observe that in \eqref{eq:state1} and \eqref{eq:costate1} there are two boundary conditions for the state $x$ while none for the co-state $\lambda$. It requires implementing specialized computational methods, such as shooting methods, to solve such a two-point boundary value problem, which in general involve intensive numerical optimizations. Here, we adopt the idea of the Sweep method
by letting
\begin{equation} 
	\label{eq:lam}
	\lambda(t) = K(t)x(t)+S(t)\nu,
\end{equation}
with $\l(t_f)=\nu$, where $K(t),~S(t)\in\mathbb{R}^{n\times n}$ and $\nu$ is the multiplier, a constant associated with the terminal constraint $\psi$, which in this case is $\psi(x(t_f))=x(t_f)=x_f$. From the transversality condition in Pontryagin's maximum principle, we know that $K(t_f)=0$ because there is no terminal cost and $S(t_f)=\frac{\partial\psi}{\partial x}\big|_{x(t_f)} = I$. Moreover, if $K$ is chosen to satisfy the Riccati equation
\begin{eqnarray}
	\label{eq:kRiccati}
	\dot{K}(t)=-Q -\tilde{A}^TK(t)-K(t)\tilde{A}+K(t)\tilde{B}R^{-1}\tilde{B}^TK(t),
\end{eqnarray}
with the terminal condition $K(t_f)=0$, then $S$ satisfies the matrix differential equation
\begin{eqnarray} \label{eq:s}
	\dot{S}(t) = -(\tilde{A}^T- K(t)\tilde{B}R^{-1}\tilde{B}^T)S(t), 
\end{eqnarray}
with the terminal condition $S(t_f)=I$, by taking the time derivative of \eqref{eq:lam} and using \eqref{eq:state1}, \eqref{eq:costate1} and \eqref{eq:kRiccati}. In addition, in order to fulfill the terminal condition $\psi(x(t_f))=x_f$ at time $t_f$, the multiplier $\nu$ associated with $\psi$ must satisfy 
\begin{eqnarray} 
	\label{eq:endpoint}
 	x_f = S^T(t)x(t)+P(t)\nu
\end{eqnarray}
for all $t\in [0,t_f]$, where $P(t)\in\mathbb{R}^{n\times n}$ obeys %is characterized by 
the matrix differential equation
\begin{align}
	\label{eq:pode}
	\dot{P}(t)-S^T(t)\tilde{B}R^{-1}\tilde{B}^TS(t)=0,
 \end{align}
with the terminal condition $P(t_f)=0$. It follows from \eqref{eq:endpoint} using $t=0$ that
\begin{equation}
	\nu=\big[P(0)\big]^{-1} \big[x_f-S^T(0)x_0\big],
\end{equation}
% which is a constant, 
provided $P(t)$ is invertible for $t\in[0,t_f]$. More details about the Sweep method based on the notion of flow mapping are provided in Appendix \ref{appd:sweep}.

% =======================================================
\subsection{Iteration procedure}
\label{sec:iterative}
The optimal solution of the problem \eqref{eq:oc1} is characterized by the homogeneous time-varying linear system described in \eqref{eq:state1} and \eqref{eq:costate1}, and we will solve for $x$ and $\l$ via an iterative procedure, which is based on analytical expressions and requires no numerical optimizations. To proceed this, we write \eqref{eq:state1} and \eqref{eq:costate1} as the iteration equations,
\begin{align}
	\label{eq:x_k}
	& \dot{x}^{(k+1)} = \tilde{A}^{(k)}x^{(k+1)} - \tilde{B}^{(k)}R^{-1}(\tilde{B}^{(k)})^T\lambda^{(k+1)},\\  
	\label{eq:lambda_k}
	& \dot{\lambda}^{(k+1)}=- \tilde{Q}^{(k)}x^{(k+1)} -(\tilde{A}^{(k)})^T\lambda^{(k+1)},
\end{align}
with identical boundary conditions $x^{(k+1)}(0)=x_0$ and $x^{(k+1)}(t_f)=x_f$ for all $k=0,1,2,\ldots$, where $\tilde{A}^{(k)}$, $\tilde{B}^{(k)}R^{-1}(\tilde{B}^{(k)})^T$, and $\tilde{Q}^{(k)}$ are defined according to %similarly as in 
\eqref{eq:A}, \eqref{eq:B}, and \eqref{eq:Q}, by
\begin{align}
	\label{eq:Ak}
	& \tilde{A}_{ij}^{(k)} = A_{ij}-\Big[(N_j R^{-1} \big(B+\sum_{j=1}^n x_j^{(k)}N_j\big)^T+\big(B+\sum_{j=1}^n x_j^{(k)}N_j) R^{-1} N_j^T\big)\lambda^{(k)}\Big]_i, \\
	\label{eq:Bk}
	&\tilde{B}^{(k)}R^{-1}(\tilde{B}^{(k)})^T = BR^{-1}B^T-\big(\sum_{j=1}^n x_j^{(k)}N_j\big)R^{-1}\big(\sum_{j=1}^n x_j^{(k)}N_j\big)^T,\\
	\label{eq:Qk}
	&\tilde{Q}^{(k)}=Q.
	% \label{eq:Ak}
	% 	& \tilde{A}^{(k)}_{ij} = A_{ij}-\Big[(N_jR^{-1}(B+\big(\sum_{j=1}^n x^{(k)}_jN_j\big))^T+(B+\big(\sum_{j=1}^n x^{(k)}_jN_j\big))R^{-1}N_j^T)\lambda^{(k)}\Big]_i, \\
	% 	\label{eq:Bk}
	% 	&\tilde{B}^{(k)}R^{-1}(\tilde{B}^{(k)})^T = BR^{-1}B^T-\big(\sum_{j=1}^n x^{(k)}_jN_j\big)R^{-1}\big(\sum_{j=1}^n x^{(k)}_jN_j\big)^T,\\
	% 	\label{eq:Qk}
	% 	&\tilde{Q}^{(k)}=Q.
\end{align}
Applying the Sweep method introduced in Section \ref{sec:sweep}, we let $\lambda^{(k+1)}(t) = K^{(k+1)}(t)x^{(k+1)}(t)+S^{(k+1)}(t)\nu^{(k+1)}$ for $t\in[0,t_f]$, where $K^{(k)}$ satisfies the Riccati equation
\begin{eqnarray}
	\label{eq:K}
	\dot{K}^{(k+1)}=-Q^{(k)}-K^{(k+1)}\tilde{A}^{(k)}-(\tilde{A}^{(k)})^TK^{(k+1)}+K^{(k+1)}\tilde{B}^{(k)}R^{-1}(\tilde{B}^{(k)})^TK^{(k+1)},
\end{eqnarray}
with the boundary condition $K^{(k+1)}(t_f)=0$, and $S^{(k)}$ follows 
\begin{align}
	\label{eq:Sk}
	\dot{S}^{(k+1)}=-\Big[(\tilde{A}^{(k)})^T-K^{(k+1)}\tilde{B}^{(k)}R^{-1}(\tilde{B}^{(k)})^T\Big] S^{(k+1)},\quad S^{(k+1)}(t_f)=I.
\end{align}
Moreover, the multiplier $\nu^{(k)}$ satisfies
\begin{equation}
	\label{eq:nuk}
	\nu^{(k+1)}=\big[P^{(k+1)}(0)\big]^{-1} \big[x_f-(S^{(k+1)})^T(0)x_0\big],
\end{equation}
where $P^{(\cdot)}(t)\in\mathbb{R}^{n\times n}$ is invertible (see Lemma \ref{lem:P_inverse} in Section \ref{sec:convergence}) and satisfies the dynamic equation
\begin{align}
	\label{eq:Pk}
	\dot{P}^{(k+1)} = (S^{(k+1)})^T\tilde{B}^{(k)}R^{-1}(\tilde{B}^{(k)})^TS^{(k+1)}
\end{align}
with the terminal condition $P^{(k+1)}(t_f)=0$. Then, the optimal control \eqref{eq:u} for the original Problem \eqref{eq:oc1} can be expressed as
\begin{equation}
	\label{eq:ocuk}
	u^*(t) = -R^{-1}\Big[B+\sum_{j=1}^n x^*_j(t)N_j\Big]^T [K^*(t)x^*(t)+S^*(t)\nu^*],
\end{equation}
if this iterative procedure is convergent, where $x^{(k)}\rightarrow x^*$, $K^{(k)}\rightarrow K^*$, and $S^{(k)}\nu^{(k)}\rightarrow S^*\nu^*$.

% ================== Remark 1 ==============
\begin{remark}
	\label{rmk:initialization}
	The iterative method can be initialized by conveniently using the optimal control of the system involving only the linear part of the bilinear system in \eqref{eq:oc1}, i.e., the LQR control. That is,
the solution $({x}^{(0)}(t),\l^{(0)}(t))$ to the homogeneous system 
\begin{align*}
	& \dot{x}^{(0)}=Ax^{(0)}-BR^{-1}B^T \lambda^{(0)}, \quad x^{(0)}(0)=x_0,\quad x^{(0)}(t_f)=x_f,\\ 
	& \dot{\lambda}^{(0)}=-A^T \lambda^{(0)}.
\end{align*}
	However, the linear system $\dot{x}=Ax+Bu$ may be uncontrollable	so that the desired transfer between $x_0$ and $x_f$ is impossible and the LQR solution does not exist. In such a case, any state trajectory with the endpoints $x_0$ and $x_f$ can be a feasible initial trajectory $x^{(0)}(t)$ of the iterative procedure.
\end{remark}

% ============================================================
\subsection{A special case: minimum-energy control of bilinear systems}
\label{eq:minimum_energy}
 Before analyzing the convergence of the iterative method, we illustrate the procedure using the example of minimum-energy control of bilinear systems, which is a special case of Problem \eqref{eq:oc1} with $Q=0$. Consider the following fixed-endpoint optimal control problem, %minimum-energy control problem involving a time-invariant bilinear system, given by
\begin{align}
	\min &\quad J=\frac{1}{2}\int_0^{t_f} u^T(t)Ru(t) \, dt, \nonumber\\
	\label{eq:ocme}
	{\rm s.t.} &\quad \dot{x}(t)=Ax+Bu+\left[\sum_{j=1}^n x_j(t) N_j\right]u, \tag{P3} \\ %\Big(\sum_{i=1}^n x_j N_j\Big) u, \\
	&\quad x(0)=x_0, \quad x(t_f)=x_f. \nonumber
\end{align}
The Hamiltonian of this %optimal control 
problem is $H(x,u,\lambda)=\frac{1}{2}u^TRu+\lambda^T[Ax+Bu+(\sum_{j=1}^n x_jN_j) u]$, where $\lambda(t)\in\mathbb{R}^n$ is the co-state vector. The optimal control is of the form as in \eqref{eq:u} and the optimal state and co-state trajectories satisfy \eqref{eq:state1} and \eqref{eq:costate1}, respectively, with $Q=0$. The respective iteration equations follow \eqref{eq:x_k} and \eqref{eq:lambda_k} with $Q^{(k)}=0$ for all $k=0,1,2,\ldots$.

Following the iterative method presented in Section \ref{sec:iterative}, we represent the costate $\lambda^{(k+1)}(t)=K^{(k+1)}(t)x^{(k+1)}(t)+S^{(k+1)}(t)\nu^{(k+1)}$, $t\in[0,t_f]$, and the matrix $K^{(k+1)}(t)\in\mathbb{R}^{n\times n}$ satisfies the Riccati equation,
\begin{align}
	\label{eq:kme}
	\dot{K}^{(k+1)} =& -K^{(k+1)}\tilde{A}^{(k)}-(\tilde{A}^{(k)})^TK^{(k+1)} + K^{(k+1)}\tilde{B}^{(k)}R^{-1}(\tilde{B}^{(k)})^TK^{(k+1)}, %\quad K^{(k+1)}(T)=0.
\end{align}
with the terminal condition $K^{(k+1)}(t_f)=0$, which has the trivial solution, $K^{(k+1)}(t)\equiv 0$, $\forall\, k=0,1,2,\ldots$, and for $t\in [0,t_f]$. This gives 
\begin{eqnarray}
	\label{eq:lambda_min_energy}
	\lambda^{(k+1)}(t)=S^{(k+1)}(t)\nu^{(k+1)},
\end{eqnarray}
and $S^{(k+1)}$ satisfies %$\dot{S}^{(k+1)} = -(\tilde{A}^{(k)})^TS^{(k+1)}$ with $S^{(k+1)}(t_f)=I$.
\begin{eqnarray}
	\label{eq:S(k+1)}
	\dot{S}^{(k+1)} = -(\tilde{A}^{(k)})^TS^{(k+1)}, \quad S^{(k+1)}(t_f)=I.
\end{eqnarray}
% which implies that $S^{(k+1)}$ is the transition matrix of the homogeneous linear system $\dot{x}=-(\tilde{A}^{(k)})^Tx$, i.e., $S^{(k+1)}(t)=\Phi_{-(\tilde{A}^{(k)})^T}(t,t_f)=\Phi^T_{\tilde{A}^{(k)}}(t_f,t)$. 
In addition, %the auxiliary variable $P^{(\cdot)}(t)\in\mathbb{R}^{n\times n}$ satisfies \eqref{eq:Pk}, and 
the multiplier associated with the terminal constraint is expressed as in \eqref{eq:nuk}. Combining \eqref{eq:lambda_min_energy} with \eqref{eq:u} gives the minimum-energy control at the $(k+1)^{th}$ iteration,
\begin{equation}
	\label{eq:ocume}
	(u^*)^{(k+1)}(t) = -R^{-1}\Big[B+\sum_{j=1}^n x^{(k+1)}_jN_j\Big]^T S^{(k+1)}\nu^{(k+1)}.
\end{equation}
Note that the auxiliary variable $P^{(k)}(t)\in\mathbb{R}^{n\times n}$ at each iteration satisfies \eqref{eq:Pk}, and thus $$P^{(k+1)}(0)=-\Phi_{\tilde{A}^{(k)}}(t_f,0) W^{(k+1)} \Phi^T_{\tilde{A}^{(k)}}(t_f,0),$$
where $\Phi^T_{\tilde{A}^{(k)}}(t_f,t)=\Phi_{-(\tilde{A}^{(k)})^T}(t,t_f)$ is the transition matrix for the homogeneous equation \eqref{eq:S(k+1)}
% $S^{(k+1)}(t)=\Phi_{-(\tilde{A}^{(k)})^T}(t,t_f)=\Phi^T_{\tilde{A}^{(k)}}(t_f,t)$. 
and
$$W^{(k+1)}=\int_0^{t_f} \Phi_{\tilde{A}^{(k)}}(0,\s) \tilde{B}^{(k)} R^{-1} (\tilde{B}^{(k)})^T \Phi^T_{\tilde{A}^{(k)}} (0,\s)d\s$$ is the controllability Gramian for the time-varying linear system as in Problem \eqref{eq:ocme}, or, equivalently, as in \eqref{eq:state1} and \eqref{eq:costate1} with $Q=0$. Moreover, the closed-loop expression in \eqref{eq:ocume} is consistent with the open-loop expression of the minimum-energy control in terms of the controllability Gramian, that is,
\begin{eqnarray}
	\label{eq:u*}
	(u^*)^{(k+1)}= R^{-1} \big[B+\sum_{j=1}^n x^{(k+1)}_jN_j\big]^T\Phi^T_{\tilde{A}^{(k)}}(0,t)\big(W^{(k+1)}\big)^{-1}\xi^{(k)},
\end{eqnarray}
where $\xi^{(k)}=\Phi_{\tilde{A}^{(k)}} (0,t_f) x_f -x_0$.

%%%%%%%%%%%%%%%%%%%%%%%%%%%%%%%%%%%%%%%%%%%%%%%%%%%%%%%%%%%%%%%%%%%%%%%%%%%
\section{Convergence of the Iterative Method}
\label{sec:convergence}
Following the iterative algorithm described in Section \ref{sec:iterative}, we expect to find the optimal control for Problem \eqref{eq:oc1}, provided the iterations are convergent. In this section, we show that the convergence of this algorithm is pertinent to the controllability of the linear system considered at each iteration and depends on the choice of the weight matrix $R$. %{\color{red} We would like to emphasize that iterative methods for solving free-endpoint optimal control problems and tracking problems involving bilinear systems were studied \cite{Hofer1988, Cimen2004} (can be moved to Introduction or removed)}. 
In Section \ref{sec:ensemble}, we will extend this iterative method to solve optimal control problems involving bilinear ensemble systems.

To facilitate the proof, we introduce the following mathematical tools. Considering the Banach spaces, $\mathcal{X}\doteq C([0,t_f];\, \mathbb{R}^n)$, $\mathcal{Y}\doteq C([0,t_f];\, \mathbb{R}^{n\times n})$, and $\mathcal{Z}\doteq C([0,t_f];\, \mathbb{R}^n)$ with the norms
\begin{align}
	\label{eq:norm1}
	\| x \|_{\alpha} &= \sup_{t\in[0,t_f]} \big[\|x(t)\|\exp(-\a t)\big], & \text{for}\ \ x\in\mathcal{X}, \\
	\label{eq:norm2}
	\| y \|_{\alpha} &= \sup_{t\in[0,t_f]} \big[\|y(t)\|\exp(-\a (t_f-t))\big], & \text{for}\ \ y\in\mathcal{Y}, \\
	\label{eq:norm3}
	\| z \|_{\a} &= \sup_{t\in[0,t_f]}\big[\|z(t)\|\exp(-\a (t_f-t))\big], & \text{for}\ \ z\in\mathcal{Z}, 
\end{align}
in which $\|v\|=\sum_{i=1}^n |v_i|$ for $v\in\mathbb{R}^n$ and $\|D\| = \max_{1\leq j\leq n} \sum_{i=1}^n |D_{ij}|$ for $D\in \mathbb{R}^{n\times n}$,
% \begin{align*}
% \|z\|=\left[\sum_{i=1}^n z_i^2 \right]^{1/2}, \qquad & z\in\mathbb{R}^n \\
% \|A\|=\left[\sum_{i,j=1}^n A_{ij}^2 \right]^{1/2}, \qquad & A\in\mathbb{R}^{n\times n} 
% \end{align*}
and the parameter $\alpha$ %that defines the norm 
serves as an additional degree of freedom to control the rate of convergence \cite{Hofer1988}, we define the operators $T_1:\mathcal{X} \times \mathcal{Y} \times \mathcal{Z} \rightarrow \mathcal{X}$, $T_2: \mathcal{X} \times \mathcal{Y} \times \mathcal{Z} \rightarrow \mathcal{Y}$, and $T_3:\mathcal{X} \times \mathcal{Y} \times \mathcal{Z} \rightarrow \mathcal{Z}$ that characterize the dynamics of $x\in\mathcal{X}$, $K\in\mathcal{Y}$, and $S\nu\in\mathcal{Z}$ % $x(t)$, $K(t)$, and $S(t)\nu$ 
as described in Section \ref{sec:iterative}, given by
\begin{align}
	\frac{d}{dt}T_1[x,K,S\nu](t) &= \tilde{A}(x(t),K(t),S(t)\nu)T_1[x,K,S\nu](t)-\tilde{B}(x(t))R^{-1}\tilde{B}^T(x(t))T_3[x,K,S\nu](t) \nonumber \\
	\label{eq:T1}
	& -\tilde{B}(x(t))R^{-1}\tilde{B}^T(x(t))T_2[x,K,S\nu](t)T_1[x,K,S\nu](t), \\
	T_1[x,K,S\nu](0) &= x_0 \nonumber \\ %, \quad T_1[x,K,S\nu](T) = x_f, \nonumber \\
	\frac{d}{dt}T_2[x,K,S\nu](t) &= - Q + T_2[x,K,S\nu](t)\tilde{B}(x(t))R^{-1}\tilde{B}^T(x(t))T_2[x,K,S\nu](t) \nonumber \\
	\label{eq:T2}
	& - T_2[x,K,S\nu](t)\tilde{A}(x(t),K(t),S(t)\nu)- \tilde{A}^T(x(t),K(t),S(t)\nu)T_2[x,K,S\nu](t) \\
	T_2[x,K,S\nu](t_f) &= 0, \nonumber \\
	\frac{d}{dt}T_3[x,K,S\nu](t) &= -\Big[ \tilde{A}^T(x(t),K(t),S(t)\nu) - T_2[x,K,S\nu](t)\tilde{B}(x(t))R^{-1}\tilde{B}^T(x(t))\Big] \cdot \nonumber \\
	\label{eq:T3}
	&\quad  T_3[x,K,S\nu](t), \\
	T_3[x,K,S\nu](t_f) &= \nu(T_1[x,K,S\nu],T_2[x,K,S\nu],T_3[x,K,S\nu]) \nonumber
\end{align}
where $\nu(T_1[x,K,S\nu],T_2[x,K,S\nu],T_3[x,K,S\nu])$ is the multiplier satisfying \eqref{eq:nuk}. With these definitions %tools 
and the following lemma, the convergence of the iterative method can be developed using the fixed-point theorem.

% ================== Lemma 3.1 ==============
\begin{lemma}
	\label{lem:P_inverse}
	The matrix $P^{(k+1)}(t)$ as in \eqref{eq:Pk} is nonsingular over $t\in [0,t_f]$ at each iteration $k$ if and only if the time-varying linear system in Problem \eqref{eq:oc2} is controllable over $[0,t_f]$ \cite{Schaettler13}.
\end{lemma}

{\it Proof:} See Appendix \ref{appd:Pnonsingular}. \hfill$\Box$

% ===================== Theorem 3.2 ===================
{\theorem
	\label{thm:convergence}
	Consider the iterative method with the iterations evolving according to 
	\begin{align}
		\label{eq:T1k}
		x^{(k+1)}(t) &= T_1[x^{(k)},K^{(k)},S^{(k)}\nu^{(k)}](t),  \\
		\label{eq:T2k}
		K^{(k+1)}(t) &= T_2[x^{(k)},K^{(k)},S^{(k)}\nu^{(k)}](t),  \\
		\label{eq:T3k}
		S^{(k+1)}(t)\nu^{(k+1)} &= T_3[x^{(k)},K^{(k)},S^{(k)}\nu^{(k)}](t),
	\end{align}
	where the operators $T_1$, $T_2$, and $T_3$ are defined in \eqref{eq:T1}, \eqref{eq:T2}, and \eqref{eq:T3}, respectively. If at each iteration $k$ the linear system as in \eqref{eq:oc2} is controllable, then $T_1$, $T_2$, and $T_3$ are contractive. Furthermore, starting with a triple of feasible trajectories $(x^{(0)},K^{(0)},S^{(0)}\nu^{(0)})$, the iteration procedure is convergent, and the sequences $x^{(k)}$, $K^{(k)}$ and $S^{(k)}\nu^{(k)}$ converge to the unique fixed points, $x^*$, $K^*$, and $(S\nu)^*$, respectively.
}
% {\theorem
% 	\label{thm:convergence}
% 	Consider the operators $T_1$, $T_2$, and $T_3$ defined in \eqref{eq:T1}, \eqref{eq:T2}, and \eqref{eq:T3}, respectively, and consider the %iteration procedure 
% 	{\color{blue} iterations evolving according to}
% 	\begin{align}
% 		\label{eq:T1k}
% 		x^{(k+1)}(t) &= T_1[x^{(k)},K^{(k)},S^{(k)}\nu^{(k)}](t),  \\
% 		\label{eq:T2k}
% 		K^{(k+1)}(t) &= T_2[x^{(k)},K^{(k)},S^{(k)}\nu^{(k)}](t),  \\
% 		\label{eq:T3k}
% 		S^{(k+1)}(t)\nu^{(k+1)} &= T_3[x^{(k)},K^{(k)},S^{(k)}\nu^{(k)}](t).
% 	\end{align}
% 	If at each iteration $k$ the linear system as in \eqref{eq:oc2} is controllable, then $T_1$, $T_2$, and $T_3$ are contractive. Furthermore, starting with a triple of feasible trajectories $(x^{(0)},K^{(0)},S^{(0)}\nu^{(0)})$, the iteration procedure defined in \eqref{eq:T1} \eqref{eq:T2} and \eqref{eq:T3} is convergent, and the sequences $x^{(k)}$, $K^{(k)}$ and $S^{(k)}\nu^{(k)}$ converge to the unique fixed points, $x^*$, $K^*$, and $(S\nu)^*$, respectively.
% }

{\it Proof:}
Because the linear system in \eqref{eq:oc2} is controllable at each iteration $k$, by Lemma \ref{lem:P_inverse} the matrix $P^{(k+1)}$ defined in \eqref{eq:Pk} is invertible and hence the multiplier $\nu^{(k+1)}$ expressed in \eqref{eq:nuk} is well-defined. Then, we have, at time $t_f$, $S^{(k+1)}(t_f)\nu^{(k+1)} = T_3[x^{(k)},K^{(k)},S^{(k)}\nu^{(k)}](t_f)=\nu^{(k+1)}$, since $S^{(k+1)}(t_f)=I$.

From \eqref{eq:Ak} and \eqref{eq:Bk}, for each fixed $t\in [0,t_f]$, we obtain the bounds
\begin{align*} 
	\|\tilde{A}^{(k+1)}-\tilde{A}^{(k)}\| &\leq \Big[\sum_{i=1}^n \|G_i \| ^2 \Big] ^{1/2} \| \lambda^{(k+1)}-\lambda^{(k)}\| \\
	& + \|\Big[\sum_{i,j=1}^n \|H_{ij} \| ^2 \Big] ^{1/2} \left\lbrace \| \lambda^{(k+1)} \| \| x^{(k+1)}-x^{(k)} \| + \| x^{(k)} \| \| \lambda^{(k+1)}-\lambda^{(k)}\| \right\}, \\
	\|\tilde{B}^{(k+1)}R^{-1}(\tilde{B}&^{(k+1)})^T-\tilde{B}^{(k)}R^{-1}(\tilde{B}^{(k)})^T\|\leq\|\Big[\sum_{i,j=1}^n \|H_{ij} \| ^2 \Big] ^{1/2} \| (x^{(k+1)})^2 - (x^{(k)})^2 \|,
\end{align*}
where $G_i = N_i R^{-1} B^T + B R^{-1} N^T_i$ and $H_{ij}=N_i R^{-1} N^T_j + N_j R^{-1} N^T_i$, and from \eqref{eq:Qk}, we have $\tilde{Q}^{(k+1)}=\tilde{Q}^{(k)}$ for all $k=0,1,2,\ldots$. Substituting \eqref{eq:lambda_k} into the above inequalities, we can write these bounds 
in terms of $\|x^{(k+1)}-x^{(k)}\|$, $\|K^{(k+1)}-K^{(k)}\|$ and $\|S^{(k+1)}\nu^{(k+1)}-S^{(k)}\nu^{(k)}\|$, given by
\begin{align} 
	\label{eq:dA}
	\|\tilde{A}^{(k+1)}-\tilde{A}^{(k)}&\| \leq \Big\{ \Big[\sum_{i=1}^n \|G_i \| ^2 \Big] ^{1/2} + \|x^{(k)} \|\Big[\sum_{i,j=1}^n \|H_{ij} \| ^2 \Big] ^{1/2} \Big\} \, \cdot \nonumber \\
	& \big\{ \| K^{(k+1)} \|\|x^{(k+1)}-x^{(k)}\| + \|K^{(k+1)}-K^{(k)}\| \|x^{(k)} \| + \|S^{(k+1)}\nu^{(k+1)}-S^{(k)}\nu^{(k)}\| \big\} \\
	&  + \Big[\sum_{i,j=1}^n \|H_{ij} \| ^2 \Big] ^{1/2} \| K^{(k+1)} \| \| x^{(k+1)} \| + \| S^{(k+1)}\nu^{(k+1)} \| \| x^{(k+1)}-x^{(k)} \|, \nonumber\\
	\label{eq:dB}
	\|\tilde{B}^{(k+1)}R^{-1} (\tilde{B} & ^{(k+1)})^T-\tilde{B}^{(k)}R^{-1}(\tilde{B}^{(k)})^T\|\leq\Big[\sum_{i,j=1}^n \|H_{ij} \|^2 \Big]^{1/2} \cdot \nonumber \\
	& \qquad\qquad\qquad\qquad\qquad\quad\Big\{\| x^{(k+1)} \| + \| x^{(k)} \| \Big\} \| x^{(k+1)}-x^{(k)} \|.
\end{align}
In addition, the solution to \eqref{eq:Sk} is given by
\begin{align} 
\nonumber
S^{(k+1)}(t) & = \Phi_{-[(\tilde{A}^{(k)})^T - K^{(k+1)}\tilde{B}^{(k)}R^{-1}(\tilde{B}^{(k)})^T]}(t,t_f)\, S^{(k+1)}(t_f) \\
\label{eq:transition}
& = \Phi^T_{[\tilde{A}^{(k)} - \tilde{B}^{(k)}R^{-1}(\tilde{B}^{(k)})^TK^{(k+1)}]}(t_f,t),
\end{align}
where $\Phi_{(.)}$ denotes the transition matrix associated with the homogeneous system \eqref{eq:Sk} and $S^{(k+1)}(t_f)=I$. Then, we have
\begin{align}
	\label{eq:dS}
 	\|S^{(k+1)}(t)-S^{(k)}(t)\| & \leq \int_t^T \|\big[ (S^{(k+1)})^T(t)\big] ^{-1}\|\|S^{(k+1)})^T(\s)\|\Big[\| \tilde{A}^{(k)}(\sigma)-\tilde{A}^{(k-1)}(\sigma) \| \nonumber\\
	& +  \| \tilde{B}^{(k)}R^{-1}(\tilde{B}^{(k)})^T(\sigma)-\tilde{B}^{(k-1)}R^{-1}(\tilde{B}^{(k-1)})^T(\sigma) \| \|K^{(k+1)}(\s)\| \\
	& + \|\tilde{B}^{(k-1)}R^{-1}(\tilde{B}^{(k-1)})^T(\sigma) \|\| K^{(k+1)}(\sigma)-K^{(k)}(\sigma) \|\Big] \|S^{(k)}(\s)\|d\sigma. \nonumber
\end{align}
From the Riccati equation for $K^{(k)}$ described in \eqref{eq:K}, we can write the differential equation for the difference $K^{(k+1)}-K^{(k)}$ as
\begin{align}
	\dfrac{d}{dt} (K^{(k+1)}- K^{(k)})& =- (K^{(k+1)}- K^{(k)})\Big[\tilde{A}^{(k)} - \tilde{B}^{(k)}R^{-1}(\tilde{B}^{(k)})^T K^{(k+1)} \Big] \nonumber\\
	\label{eq:Kd}
	& - \Big[\tilde{A}^{(k)} - \tilde{B}^{(k)}R^{-1}(\tilde{B}^{(k)})^T K^{(k+1)} \Big]^T (K^{(k+1)}- K^{(k)}) \\
	& - K^{(k)}(\tilde{A}^{(k)}-\tilde{A}^{(k-1)}) - (\tilde{A}^{(k)}-\tilde{A}^{(k-1)})^T K^{(k+1)} \nonumber \\
	& + K^{(k)}( \tilde{B}^{(k)}R^{-1}(\tilde{B}^{(k)})^T-\tilde{B}^{(k-1)}R^{-1}(\tilde{B}^{(k-1)})^T ) K^{(k+1)}, \nonumber
\end{align}
with the terminal condition $K^{(k+1)}(t_f)-K^{(k)}(t_f)=0$. Applying the variation of constants formula, backward in time from $t=t_f$, to \eqref{eq:Kd} and employing \eqref{eq:transition} yield
\begin{align*}
	K^{(k+1)}(t)- & K^{(k)}(t) = (S^{(k)})^T(t)\Big\{ \int_t^{t_f} \Big[(S^{(k)})^T(\s)\Big]^{-1} \Big[K^{(k)}(\s)\big(\tilde{A}^{(k)}(\s)-\tilde{A}^{(k-1)}(\s)\big) \\
	& +\big(\tilde{A}^{(k)}(\s)-\tilde{A}^{(k-1)}(\s)\big)^T K^{(k+1)}(\s)-K^{(k)}(\s)\cdot \\
	& \big(\tilde{B}^{(k)}R^{-1}(\tilde{B}^{(k)})^T-\tilde{B}^{(k-1)}R^{-1}(\tilde{B}^{(k-1)})^T\big) K^{(k+1)}(\s)\Big] \Big[S^{(k+1)}(\s)\Big]^{-1} d\s \Big\}  S^{(k+1)}(t),
\end{align*}
which results in
\begin{align} 
	\label{eq:dK}
	\|K^{(k+1)}(t)-K^{(k)}(t)\| & \leq \int_t^{t_f} \Big[\beta_1 \|\tilde{A}^{(k)}(\s)-\tilde{A}^{(k-1)}(\s)\| \nonumber \\
& +\beta_2 \|\tilde{B}^{(k)}R^{-1}(\tilde{B}^{(k)})^T(\s)-\tilde{B}^{(k-1)}R^{-1}(\tilde{B}^{(k-1)})^T(\s)\| \Big] d\s,
\end{align}
where $\beta_1$ and $\beta_2$ are both finite time-varying coefficients (see Appendix \ref{appd:betas}).

Similarly, from \eqref{eq:x_k} and \eqref{eq:lambda_k}, we can write the differential equation for $(x^{(k+1)}-x^{(k)})$, that is,
\begin{align}
	\frac{d}{dt} (x^{(k+1)}- x^{(k)}) &= \Big[\tilde{A}^{(k)} - \tilde{B}^{(k)}R^{-1}(\tilde{B}^{(k)})^TK^{(k+1)} \Big] (x^{(k+1)}-x^{(k)}) \nonumber\\
	&+ \Big\{ (\tilde{A}^{(k)}-\tilde{A}^{(k-1)}) - \left( \tilde{B}^{(k)}R^{-1}(\tilde{B}^{(k)})^T-\tilde{B}^{(k-1)}R^{-1}(\tilde{B}^{(k-1)})^T \right) K^{(k+1)}  \nonumber\\
	\label{eq:xd}
	& - \tilde{B}^{(k-1)}R^{-1}(\tilde{B}^{(k-1)})^T (K^{(k+1)}-K^{(k)})\Big\} x^{(k)} \\
	& - \left( \tilde{B}^{(k)}R^{-1}(\tilde{B}^{(k)})^T-\tilde{B}^{(k-1)}R^{-1}(\tilde{B}^{(k-1)})^T \right) S^{(k+1)}\nu^{(k+1)} \nonumber \\
	& - \tilde{B}^{(k-1)}R^{-1}(\tilde{B}^{(k-1)})^T \big(S^{(k+1)}\nu^{(k+1)}-S^{(k)}\nu^{(k)}\big), \nonumber 
\end{align}
with the terminal condition $x^{(k+1)}(t_f)-x^{(k)}(t_f)=0$. Applying the variation of constants formula to \eqref{eq:xd} yields,
\begin{align*}
	x^{(k+1)}(t)-x^{(k)}(t) &= \Big[ (S^{(k+1)})^T(t)\Big]^{-1}\int_0^t (S^{(k+1)})^T(\s)\Big\{ \Big[\big(\tilde{A}^{(k)}(\s)-\tilde{A}^{(k-1)}(\s)\big) \\
	& - \big(\tilde{B}^{(k)}R^{-1}(\tilde{B}^{(k)})^T(\sigma)-\tilde{B}^{(k-1)}R^{-1}(\tilde{B}^{(k-1)})^T(\sigma)\big) K^{(k+1)}(\sigma) \\
	& - \tilde{B}^{(k-1)}R^{-1}(\tilde{B}^{(k-1)})^T (\sigma) \big(K^{(k+1)}(\sigma)-K^{(k)}(\sigma)\big) \Big] x^{(k)}(\sigma) \\
	& - \left( \tilde{B}^{(k)}R^{-1}(\tilde{B}^{(k)})^T-\tilde{B}^{(k-1)}R^{-1}(\tilde{B}^{(k-1)})^T \right) S^{(k+1)}\nu^{(k+1)} \nonumber \\
	& - \tilde{B}^{(k-1)}R^{-1}(\tilde{B}^{(k-1)})^T(S^{(k+1)}\nu^{(k+1)}-S^{(k)}\nu^{(k)}) \Big\} d\s.
\end{align*}
It follows that
\begin{align}
	\| x^{(k+1)}(t)-x^{(k)}(t) \| & \leq \int_0^t \Big[ \beta_3 \| \tilde{A}^{(k)}(\sigma)-\tilde{A}^{(k-1)}(\sigma) \| + \beta_4 \| K^{(k+1)}(\sigma) - K^{(k)}(\sigma) \| \nonumber \\
	\label{eq:dx}
	& + \beta_5 \|\tilde{B}^{(k)}R^{-1}(\tilde{B}^{(k)})^T(\sigma)-\tilde{B}^{(k-1)}R^{-1}(\tilde{B}^{(k-1)})^T (\s) \| \\
	& + \beta_6 \|S^{(k+1)}\nu^{(k+1)}(\sigma)-S^{(k)}\nu^{(k)}(\sigma)\| \Big] d\sigma, \nonumber
\end{align}
where $\beta_3$, $\beta_4$, $\beta_5$ and $\beta_6$ are all finite time-varying coefficients (see Appendix \ref{appd:betas}). Furthermore, since $\nu^{(k+1)}$ is a constant within each iteration $k$, from \eqref{eq:Sk} we can write 
$$\dfrac{d}{dt}(S^{(k+1)}\nu^{(k+1)}) = -\Big[(\tilde{A}^{(k)})^T - K^{(k+1)}\tilde{B}^{(k)}R^{-1}(\tilde{B}^{(k)})^T\Big] S^{(k+1)}\nu^{(k+1)},$$ 
with the terminal condition $S^{(k+1)}(t_f)\nu^{(k+1)} = \nu^{(k+1)}$. This allows us to write %Hence, the differential equation of the difference $(S^{(k+1)}\nu^{(k+1)}-S^{(k)}\nu^{(k)})$ is given as,
\begin{align*}
	\dfrac{d}{dt}(S^{(k+1)}\nu^{(k+1)}-S^{(k)}\nu^{(k)}) &= -\Big[\tilde{A}^{(k)} - \tilde{B}^{(k)}R^{-1}(\tilde{B}^{(k)})^TK^{(k+1)} \Big]^T (S^{(k+1)}\nu^{(k+1)}-S^{(k)}\nu^{(k)})\\
	& -\Big\{ (\tilde{A}^{(k)}-\tilde{A}^{(k-1)}) - \tilde{B}^{(k-1)}R^{-1}(\tilde{B}^{(k-1)})^T (K^{(k+1)}-K^{(k)})\nonumber\\
	& -\left( \tilde{B}^{(k)}R^{-1}(\tilde{B}^{(k)})^T-\tilde{B}^{(k-1)}R^{-1}(\tilde{B}^{(k-1)})^T \right) K^{(k+1)}  \Big\} ^T S^{(k)}\nu^{(k)},
\end{align*}
with the terminal condition $S^{(k+1)}(t_f)\nu^{(k+1)}-S^{(k)}(t_f)\nu^{(k)} = \nu^{(k+1)}-\nu^{(k)}$, and then %which gives the solution
\begin{align*}
	S^{(k+1)}(t)\nu^{(k+1)} &- S^{(k)}(t)\nu^{(k)} = \Big[ (S^{(k+1)})^T(t)\Big]^{-1} \Big\{(\nu^{(k+1)}-\nu^{(k)}) \\
	&- \int_t^{t_f} (S^{(k+1)})^T(\s) \Big[(\tilde{A}^{(k)}-\tilde{A}^{(k-1)}) - \tilde{B}^{(k-1)}R^{-1}(\tilde{B}^{(k-1)})^T (K^{(k+1)}-K^{(k)})\nonumber\\
	&- \left( \tilde{B}^{(k)}R^{-1}(\tilde{B}^{(k)})^T-\tilde{B}^{(k-1)}R^{-1}(\tilde{B}^{(k-1)})^T \right) K^{(k+1)}\Big]^T S^{(k)}\nu^{(k)} d\sigma \Big\}.
\end{align*}
From \eqref{eq:nuk} we may obtain
\begin{align}
	\label{eq:dnu}
	\nonumber
	\|\nu^{(k+1)}-\nu^{(k)}\| & \leq \|(P^{(k+1)})^{-1}\| \, \|P^{(k+1)}(t)-P^{(k)}(t)\| \, \|\nu^{(k)}\| \\ 
	& + \|(P^{(k+1)})^{-1}\| \, \|S'^{(k+1)}(t)-S'^{(k)}(t)\| \, \|x_0\|,
\end{align}
in which, by evolving \eqref{eq:Pk} backward in time from $t=t_f$, the difference $P^{(k+1)}(t)-P^{(k)}(t)$ satisfies
\begin{align}
	\nonumber
	P^{(k+1)}(t)-P^{(k)}(t) &= -\int_t^{t_f} \bigg[\big(S^{(k+1)}+S^{(k)}\big) \tilde{B}^{(k-1)}R^{-1}(\tilde{B}^{(k-1)})^T \big(S^{(k+1)}(\s)-S^{(k)}(\s)\big) \\
\label{eq:dP}
& + S^{(k+1)}\Big(\tilde{B}^{(k)}R^{-1}(\tilde{B}^{(k)}(\s))^T-\tilde{B}^{(k-1)}R^{-1}(\tilde{B}^{(k-1)}(\s))^T\Big)S^{(k)} \bigg] d\s.
\end{align}
Using \eqref{eq:dnu} and \eqref{eq:dP}, we obtain
\begin{align}
	\| S^{(k+1)}(t)\nu^{(k+1)}-S^{(k)}(t)\nu^{(k)} \| & \leq \int_t^{t_f} \Big[\beta_7 \, \| \tilde{A}^{(k)}(\sigma)-\tilde{A}^{(k-1)}(\sigma) \| \nonumber\\
	\label{eq:dcostate1}
	& + \b_8 \, \| \tilde{B}^{(k)}R^{-1}(\tilde{B}^{(k)})^T(\sigma)-\tilde{B}^{(k-1)}R^{-1}(\tilde{B}^{(k-1)})^T(\sigma) \| \\
	& + \b_9 \, \| K^{(k+1)}(\sigma)-K^{(k)}(\sigma) \| \Big] d\sigma, \nonumber
\end{align}
%\begin{align}
%	& \| S^{(k+1)}(t)\nu^{(k+1)}-S^{(k)}(t)\nu^{(k)} \| \leq \int_t^{t_f} \Big[\beta_7 \, \| \tilde{A}^{(k)}(\sigma)-\tilde{A}^{(k-1)}(\sigma) \| \nonumber\\
%	\label{eq:dcostate1}
%	& + \b_8 \, \| \tilde{B}^{(k)}R^{-1}(\tilde{B}^{(k)})^T(\sigma)-\tilde{B}^{(k-1)}R^{-1}(\tilde{B}^{(k-1)})^T(\sigma) \| + \b_9 \, \| K^{(k+1)}(\sigma)-K^{(k)}(\sigma) \| \Big] d\sigma, %\nonumber
%\end{align}

where $\beta_7$, $\beta_8$ and $\beta_9$ are finite time-varying coefficients (see Appendix \ref{appd:betas}).

% Evolving \eqref{eq:Pk} backward in time from $t=t_f$ yields
% \begin{align}
% 	\nonumber
% 	P^{(k+1)}(t)-P^{(k)}(t) &= -\int_t^{t_f} \bigg[\big(S^{(k+1)}+S^{(k)}\big) \tilde{B}^{(k-1)}R^{-1}(\tilde{B}^{(k-1)})^T \big(S^{(k+1)}(\s)-S^{(k)}(\s)\big) \\
% \label{eq:dP}
% & + S^{(k+1)}\Big(\tilde{B}^{(k)}R^{-1}(\tilde{B}^{(k)}(\s))^T-\tilde{B}^{(k-1)}R^{-1}(\tilde{B}^{(k-1)}(\s))^T\Big)S^{(k)} \bigg] d\s.
% \end{align}
% Also, it follows from \eqref{eq:nuk} that 
% \begin{align}
% 	\label{eq:dnu}
% 	\nonumber
% 	\|\nu^{(k+1)}-\nu^{(k)}\| & \leq \|(P^{(k+1)})^{-1}\| \, \|P^{(k+1)}(t)-P^{(k)}(t)\| \, \|\nu^{(k)}\| \\ 
% 	& + \|(P^{(k+1)})^{-1}\| \, \|S'^{(k+1)}(t)-S'^{(k)}(t)\| \, \|x_0\|.
% \end{align} 
% Using \eqref{eq:dP} and \eqref{eq:dnu}, we obtain %Thus, obtaining an estimate using \eqref{eq:dS}, \eqref{eq:dP} and \eqref{eq:dnu},
% \begin{align}
% 	\| S^{(k+1)}(t)\nu^{(k+1)}-S^{(k)}(t)\nu^{(k)} \| & \leq \int_t^{t_f} \Big[\beta_7 \, \| \tilde{A}^{(k)}(\sigma)-\tilde{A}^{(k-1)}(\sigma) \| \nonumber\\
% 	\label{eq:dcostate1}
% 	& + \b_8 \, \| \tilde{B}^{(k)}R^{-1}(\tilde{B}^{(k)})^T(\sigma)-\tilde{B}^{(k-1)}R^{-1}(\tilde{B}^{(k-1)})^T(\sigma) \| \\
% 	& + \b_9 \, \| K^{(k+1)}(\sigma)-K^{(k)}(\sigma) \| \Big] d\sigma, \nonumber
% \end{align}
% where $\beta_7$, $\beta_8$ and $\beta_9$ are finite time-varying coefficients (see Appendix \ref{appd:betas}).

Combining the bounds in \eqref{eq:dA}, \eqref{eq:dB}, \eqref{eq:dK}, \eqref{eq:dx}, and \eqref{eq:dcostate1}, and using the definitions of the operators $T_1$, $T_2$, and $T_3$ in \eqref{eq:T1}, \eqref{eq:T2} and \eqref{eq:T3}, respectively, we reach the inequality that holds component-wise, given by
\begin{align}
	\label{eq:contraction}
	& \left[\begin{array}{c}\| T_1[x^{(k)},K^{(k)},S^{(k)}\nu^{(k)}]-T_1[x^{(k-1)},K^{(k-1)},S^{(k-1)}\nu^{(k-1)}]\|_{\a} \\
	\| T_2[x^{(k)},K^{(k)},S^{(k)}\nu^{(k)}]-T_2[x^{(k-1)},K^{(k-1)},S^{(k-1)}\nu^{(k-1)}]\|_{\a} \\
	\|T_3[x^{(k)},K^{(k)},S^{(k)}\nu^{(k)}]-T_3[x^{(k-1)},K^{(k-1)},S^{(k-1)}\nu^{(k-1)}]\|_{\a} \end{array}\right] \nonumber \\
	& \leq M \left[\begin{array}{c}\|x^{(k)}-x^{(k-1)}\|_{\a} \\
	\|K^{(k)}-K^{(k-1)}\|_{\a} \\
	\|S^{(k)}\nu^{(k)}-S^{(k-1)}\nu^{(k-1)}\|_{\a} \end{array}\right],
\end{align}
where $M\in\mathbb{R}^{3\times 3}$ whose elements are all related to $R^{-1}$ (see Appendix \ref{appd:M}). As a result, the eigenvalues of $M$ can be made  within the unit circle by choosing sufficiently large $R$. Therefore, the operators $T_1$, $T_2$ and $T_3$ are contractive, and the fixed-point theorem \cite{Brockett70} warrants the convergence of the iterative procedure to the unique fixed points, i.e., $x^{(k)}\rightarrow x^*$, $K^{(k)}\rightarrow K^*$, and $S^{(k)}\nu^{(k)}\rightarrow (S\nu)^*$. %, which also lead to $\lambda^{(k)}\rightarrow \lambda^*$ from \eqref{eq:lam} and {\color{blue} by} the continuity of all the variables involved. 
Note that the choice of $R$ determines the magnitude of eigenvalues of $M$ and thus can also be used to improve the convergence rate of the iterative procedure. \hfill$\Box$

% As a result, the eigenvalues of $M$ can be made arbitrarily small, in particular within the unit circle, by choosing sufficiently large $R$. Therefore, the operators $T_1$, $T_2$ and $T_3$ are contractive, and the Fixed-Point theorem \cite{Brockett70} warrants the convergence of the presented iterative procedure to the unique fixed points, i.e., $x^{(k)}\rightarrow x^*$, $K^{(k)}\rightarrow K^*$, and $S^{(k)}\nu^{(k)}\rightarrow (S\nu)^*$, which also lead to $\lambda^{(k)}\rightarrow \lambda^*$ from \eqref{eq:lam} and the continuity of all the variables involved. \hfill$\Box$

% ================ Remark 2 =============
\begin{remark}[Optimality of the Convergent Solution]
	\label{rmk:necessary}
	The convergence of $x^{(k)}\rightarrow x^*$, $K^{(k)}\rightarrow K^*$, and $S^{(k)}\nu^{(k)}\rightarrow (S\nu)^*$ immediately leads to $\lambda^{(k)}\rightarrow \lambda^*$ by \eqref{eq:lam} and by the continuity of all the variables involved. This in turn guarantees that the fixed points $x^*$ and $\lambda^*$ resulting from the iterative procedure are the solutions to \eqref{eq:state1} and \eqref{eq:costate1} with the convergent $\tilde{A}$ and $\tilde{B}R^{-1}\tilde{B}^T$, denoted $A^*$ and $B^*R^{-1}(B^*)^T$, obtained from \eqref{eq:A} and \eqref{eq:B}, %\eqref{eq:state} and \eqref{eq:costate}, 
	respectively. This implies that the convergent solution pair $(x^*,\lambda^*)$ satisfies the necessary optimality condition, and thus the convergent control $u^*$ is a candidate of the optimal control for Problem \eqref{eq:oc1}.
\end{remark}

%%%%%%%%%%%%%%%%%%%%%%%%%%%%%%%%%%%%%%%%%%%%%%%%%%%%%%%%%%%%%%%%%%%%%%%%%%%
\section{Global Optimality of the Convergent Solution}
\label{sec:optimality}
We have shown in Remark \ref{rmk:necessary} that the convergent optimal control $u^*$ generated by the iterative procedure satisfies the necessary optimality condition. In this section, we will further illustrate that $u^*$ may be the unique global optimal control given appropriate assumptions on the value function associated with Problem \eqref{eq:oc2}.

% =========================================
\subsection{Optimality of the solution at each iteration}
\label{sec:sufficient}
For each iteration $k$, \eqref{eq:oc2} is a time-dependent problem with a specified the time horizon, and the optimal control satisfies the Hamilton-Jacobi-Bellman (HJB) equation, given by
\begin{align} 
	\label{eq:HJBlqr}
	\dfrac{\partial V^{(k)}}{\partial t}(t,x^{(k)}) + \min_{u\in \mathcal{U}} \Big\{ \dfrac{\partial V^{(k)}}{\partial x^{(k)}}(t,x^{(k)})^T (\tilde{A}^{(k-1)}x^{(k)}+\tilde{B}^{(k-1)}u)+\dfrac{1}{2} [(x^{(k)})^TQx^{(k)}+u^TRu]\Big\}\equiv 0,
\end{align}
with the boundary condition $V^{(k)}(t_f,x^{(k)})=0$, where $V^{(k)}$ is the value function and $\mathcal{U}$ is the set of all admissible controls. % defined over some interval $[t,T]$
Since the matrix $R\in \mathbb{R}^{n\times n}$ is positive definite, the function to be minimized in \eqref{eq:HJBlqr} is strictly convex in the control variable $u$. As a result, the minimization problem in \eqref{eq:HJBlqr} has a unique solution given by the stationary point, satisfying $(\tilde{B}^{(k-1)})^T \dfrac{\partial V^{(k)}}{\partial x^{(k)}}(t,x^{(k)}) + Ru_*^{(k)}(t)=0$, or, equivalently,
\begin{align}
	\label{eq:optimalu}
	u_*^{(k)}(t) = -R^{-1}(\tilde{B}^{(k-1)})^T \dfrac{\partial V^{(k)}}{\partial x^{(k)}}(t,x^{(k)}).
\end{align}
Substituting \eqref{eq:optimalu} into %the HJB equation in 
\eqref{eq:HJBlqr} gives a first-order nonlinear partial differential equation,
\begin{align}
	\dfrac{\partial V^{(k)}}{\partial t}(t,x^{(k)}) &+ \dfrac{1}{2} \Big[ \dfrac{\partial V^{(k)}}{\partial x^{(k)}}(t,x^{(k)})^T\tilde{A}^{(k-1)}x^{(k)} +(x^{(k)})^T(\tilde{A}^{(k-1)})^T\dfrac{\partial V^{(k)}}{\partial x^{(k)}}(t,x^{(k)})\Big] \nonumber \\
	\label{eq:fopde}
& + \dfrac{1}{2} (x^{(k)})^TQx^{(k)}-\dfrac{1}{2} \dfrac{\partial V^{(k)}}{\partial x^{(k)}}(t,x^{(k)})^T \tilde{B}^{(k-1)}R^{-1}(\tilde{B}^{(k-1)})^T \dfrac{\partial V^{(k)}}{\partial x^{(k)}}(t,x^{(k)})\equiv 0.
\end{align}
% \begin{align}
% 	\dfrac{\partial V^{(k)}}{\partial t}(t,x^{(k)}) &+ \dfrac{1}{2} \Big[ \dfrac{\partial V^{(k)}}{\partial x^{(k)}}(t,x^{(k)})^T\tilde{A}^{(k-1)}x^{(k)} +(x^{(k)})^T(\tilde{A}^{(k-1)})^T\dfrac{\partial V^{(k)}}{\partial x^{(k)}}(t,x^{(k)})\Big] + \dfrac{1}{2} (x^{(k)})^TQx^{(k)} \nonumber \\
% 	\label{eq:fopde}
% & -\dfrac{1}{2} \dfrac{\partial V^{(k)}}{\partial x^{(k)}}(t,x^{(k)})^T \tilde{B}^{(k-1)}R^{-1}(\tilde{B}^{(k-1)})^T \dfrac{\partial V^{(k)}}{\partial x^{(k)}}(t,x^{(k)})\equiv 0.
% \end{align}
Due to the quadratic and symmetric nature of \eqref{eq:fopde}, we consider the value function of the form,
\begin{align}
	\label{eq:Vfunc}
	V^{(k)}(t,x) = \dfrac{1}{2}x^TK^{(k)}x + x^TS^{(k)}\nu^{(k)} +\dfrac{1}{2}(\nu^{(k)})^TP^{(k)}\nu^{(k)},
\end{align}
with the boundary condition $V^{(k)}(t_f,x^{(k)})=0$, where $K^{(k)}(t),~S^{(k)}(t),~P^{(k)}(t)\in \mathbb{R}^{n\times n}, \forall~ t\in [0,t_f]$ and $\nu^{(k)}\in\mathbb{R}$. It is straightforward to verify that $(V^{(k)}(t,x^{(k)}),u_*^{(k)})$ %, where $V^{(k)}$ as in \eqref{eq:Vfunc} and $u_*^{(k)}$ as in \eqref{eq:optimalu}, 
is a classical solution to the HJB equation \eqref{eq:fopde} if the matrices $K^{(k)}(t)$, $S^{(k)}(t)$, and $P^{(k)}(t)$ satisfy the matrix differential equations \eqref{eq:kRiccati}, \eqref{eq:s} and \eqref{eq:pode}, respectively, with the respective boundary conditions $K^{(k)}(t_f)=0$, $S^{(k)}(t_f)=I$, and $P^{(k)}(t_f)=0$. Note that the value function $V^{(k)}(t,x^{(k)})$ in \eqref{eq:Vfunc} is continuously differentiable on its domain and extends continuously onto the terminal manifold $\mathcal{N}^{(k)} = \{ x^{(k)}(t_f) = x_f \}$.

% ==================================================================
\subsection{Global optimality of the convergent solution}
\label{sec:global}
We showed in Section \ref{sec:sufficient} that the control $u_*^{(k)}$ presented in \eqref{eq:optimalu} is the global optimum for the $k^{th}$ iteration in Problem \eqref{eq:oc2}. Here, we will show that the convergent solution $u^*$ of $u_*^{(k)}$, i.e., $u_*^{(k)}\rightarrow u^*$, %resulting from} the iterative procedure applied to solve 
for the linear problem \eqref{eq:oc2} is indeed the global optimal control for the original bilinear problem \eqref{eq:oc1} under some regularity conditions on the value function associated with \eqref{eq:oc2} expressed in \eqref{eq:Vfunc}.

% =============== Theorem 4.1 ============
\begin{theorem}
	\label{thm:globalsufficient}
	Consider the iterative method applied to Problem \eqref{eq:oc2}, and suppose that at each iteration $k$ the linear system as in \eqref{eq:oc2} is controllable. Let $u^*$ be the convergent solution of the optimal control sequence $\{u_*^{(k)}\}$ generated by the iterative procedure for $k\in\mathbb{N}$, i.e., $u_*^{(k)}\rightarrow u^*$, and let $V^*$ be the corresponding convergent value function defined in \eqref{eq:Vfunc}, i.e., $V^{(k)} \rightarrow V^*$. If (i) $V^*\in C^1$, and $\dfrac{\partial V^*}{\partial t}$ and $\dfrac{\partial V^*}{\partial x}$ are Lipschitz continuous; and (ii) there exist real-valued $L_1$ functions, $g(t)$ and $h_i(x)$, $i = 1,2,\ldots,n$, i.e., $g \in L_1([0,t_f])$ and $h_i \in L_1(M)$ where $M \subset \mathbb{R}^n$, such that $\Big|\dfrac{\partial V^{(k)}}{\partial t}(t,x^{(k)}(t))\Big| \leq g(t)$ for all $k\in\mathbb{N}$ and $t\in [0,t_f]$, and for each component $i$, $\Big| \Big[\dfrac{\partial V^{(k)}}{\partial x}(t,x)\Big]_i\Big| \leq h_i(x)$ for all $k\in\mathbb{N}$ and for all $x=x^{(k)}(t)\in\mathbb{R}^n$, then $u^*$ is a global optimum for the original Problem \eqref{eq:oc1}.
	% such that $\Big|\dfrac{\partial V^{(k)}}{\partial t}(t,x^{(k)}(t))\Big| \leq g(t)$ for all $k\in\mathbb{N}$ and $t\in [0,t_f]$, and for each component $i$, $\Big| \Big[\dfrac{\partial V^{(k)}}{\partial x}(t,x)\Big]_i\Big| \leq h_i(x)$ for all $k\in\mathbb{N}$ and for all $x=x^{(k)}(t)\in\mathbb{R}^n$, then $u^*$ is a global optimum for the original Problem \eqref{eq:oc1}.
\end{theorem}

{\it Proof:}
First of all, the conditions in (i) guarantee the existence of the optimal control. Because the linear system as in \eqref{eq:oc2} is controllable at each iteration $k$, there exist unique fixed points for the sequences $x^{(k)}$, $K^{(k)}$, and $S^{(k)}\nu^{(k)}$ such that $x^{(k)}\rightarrow x^*$, $K^{(k)}\rightarrow K^*$ and $S^{(k)}\nu^{(k)}\rightarrow (S\nu)^*$ by Theorem \ref{thm:convergence}. It follows that the partial derivatives of $V^{(k)}$ with respect to $t$ and $x$ are convergent, denoted $\dfrac{\partial V^{(k)}}{\partial t}(t,x) \rightarrow V_t(t,x)$ and $\dfrac{\partial V^{(k)}}{\partial x}(t,x)\rightarrow V_x(t,x)$, where
\begin{align*}
	V_t(t,x^*) &= \frac{1}{2}(x^*)^T \Big[ - Q - (\tilde{A}^*)^TK^* - K^*\tilde{A}^* + K^*\tilde{B}^* R^{-1}(\tilde{B}^*)^TK^*\Big] x^* \\
	&+ (x^*)^T \Big[ -(\tilde{A}^*)^T - K^*\tilde{B}^* R^{-1}(\tilde{B}^*)^T \Big]S^* \nu^* + \dfrac{1}{2} (\nu^*)^T (S^*)^T \tilde{B}^* R^{-1}(\tilde{B}^*)^TS^* \nu^*,\\
	V_x(t,x^*) &= K^*x^* +S^* \nu^*= \lambda^*,
\end{align*}
in which $\tilde{A}^*$ and $\tilde{B}^*R^{-1}(\tilde{B}^*)^T$ are the limits of $\tilde{A}^{(k)}$ and $\tilde{B}^{(k)} R^{-1}(\tilde{B}^{(k)})^T$, respectively, %$\tilde{A}^{(k)}\rightarrow\tilde{A}^*$ and $\tilde{B}^{(k)} R^{-1}(\tilde{B}^{(k)})^T\rightarrow\tilde{B}^*R^{-1}(\tilde{B}^*)^T$ 
following \eqref{eq:Ak} and \eqref{eq:Bk}. Because $\dfrac{\partial V^{(k)}}{\partial t}(t,x^{(k)}(t))$ and $\dfrac{\partial V^{(k)}}{\partial x}(t,x^{(k)}(t))$ are dominated by $g(t)$ and $h(x)=(h_1,\ldots,h_n)'$, respectively, by the Lebesgue Dominated Convergence theorem, we have
\begin{align}
	\label{eq:V*1}
	\lim_{k\rightarrow \infty} V^{(k)}(t,x) &= \lim_{k\rightarrow \infty}\int_0^t \dfrac{\partial V^{(k)}}{\partial \s}(\s,x) d\s =\int_0^t \lim_{k\rightarrow \infty}\dfrac{\partial V^{(k)}}{\partial \s}(\s,x) d\s = \int_0^t V_\s (\s,x) d\s, \\
	\label{eq:V*2}
	% \lim_{k\rightarrow \infty} V^{(k)}(t,x) 
	\lim_{k\rightarrow \infty} V^{(k)}(t,x) &= \lim_{k\rightarrow \infty}\int_{x(0)}^{x(t)} \dfrac{\partial V^{(k)}}{\partial x}(t,x) dx =\int_{x(0)}^{x(t)} \lim_{k\rightarrow \infty}\dfrac{\partial V^{(k)}}{\partial x}(t,x) dx = \int_{x(0)}^{x(t)} V_x (t,x) dx.
\end{align}
where $\lim_{k\rightarrow \infty} V^{(k)}(t,x)=V^* (t,x)$ by assumption.
% \begin{align}
% 	\label{eq:V*1}
% 	V^* (t,x) & = \lim_{k\rightarrow \infty} V^{(k)}(t,x) = \lim_{k\rightarrow \infty}\int_0^t \dfrac{\partial V^{(k)}}{\partial \s}(\s,x) d\s =\int_0^t \lim_{k\rightarrow \infty}\dfrac{\partial V^{(k)}}{\partial \s}(\s,x) d\s = \int_0^t V_\s (\s,x) d\s, \\
% 	\label{eq:V*2}
% 	V^* (t,x) & = \lim_{k\rightarrow \infty} V^{(k)}(t,x) = \lim_{k\rightarrow \infty}\int_{x(0)}^{x(t)} \dfrac{\partial V^{(k)}}{\partial x}(t,x) dx =\int_{x(0)}^{x(t)} \lim_{k\rightarrow \infty}\dfrac{\partial V^{(k)}}{\partial x}(t,x) dx = \int_{x(0)}^{x(t)} V_x (t,x) dx.
% \end{align}
Because $V^*$ is continuously differentiable with respect to both $t$ and $x$, we obtain from \eqref{eq:V*1} and \eqref{eq:V*2} the partial derivatives
\begin{align} 
	\label{eq:direvative}
	\dfrac{\partial V^*}{\partial t}(t,x) = V_t (t,x),\quad \dfrac{\partial V^*}{\partial x}(t,x) = V_x (t,x).
\end{align}
In addition, due to the convergence of the iterative procedure, \eqref{eq:fopde} is convergent to
\begin{align*}
	V_t(t,x^*) &+ V_x(t,x^*)^T\tilde{A}^*x^* + \dfrac{1}{2} (x^*)^TQx^* - \frac{1}{2} V_x(t,x^*)^T \tilde{B}^*R^{-1}(\tilde{B}^*)^T V_x(t,x^*)\equiv 0.
\end{align*}
which, by employing \eqref{eq:direvative} and $V_x(t,x^*) = \lambda^*$, can be rewritten as 
\begin{align*}
	\dfrac{\partial V^*}{\partial t}(t,x^*) +  \dfrac{\partial V^*}{\partial x}(t,x^*)^T (\tilde{A}^*x^*-\tilde{B}^*R^{-1}(\tilde{B}^*)^T\lambda^*)+\dfrac{1}{2}\big[(x^*)^TQx^*+(\lambda^*)^T \tilde{B}^*R^{-1}(\tilde{B}^*)^T\lambda^*\big]\equiv 0,
\end{align*}
with the boundary condition $V^*(t_f,x^*)=0$. Because the convergent solution pair $(x^*,\lambda^*)$ satisfies the necessary condition (see Remark \ref{rmk:necessary}), the above equation is equivalent to, by \eqref{eq:optimalu}, \eqref{eq:state1}, \eqref{eq:costate1}, and \eqref{eq:linear},
\begin{align} 
	\label{eq:HJBu}
	\dfrac{\partial V^*}{\partial t}(t,x^*) +  \dfrac{\partial V^*}{\partial x}(t,x^*)^T \Big[Ax^*+\Big(B+\sum_{i=1}^n x^*_i N_i\Big) u^*\Big]+\dfrac{1}{2} [(x^*)^TQx^*+(u^*)^TRu^* ]\equiv 0.
\end{align}
Since $V^*$ is differentiable, according to the dynamic programming principle \cite{Schaettler13}, the quantity on the left-hand side in \eqref{eq:HJBu} is non-negative for every control $u$ in the admissible control set $\mathcal{U}\subset\mathbb{R}^m$. It follows that $(V^*,u^*)$ is a solution to the HJB equation of the original Problem \eqref{eq:oc1}, that is,
\begin{align} 
	\label{eq:HJBbilinear}
	\dfrac{\partial V}{\partial t}(t,x) + \min_{u\in\mathcal{U}} \Big\{ \dfrac{\partial V}{\partial x}(t,x)^T \Big[Ax+\Big(B+\sum_{i=1}^n x_i N_i\Big) u\Big]+\dfrac{1}{2} [x^TQx+u^TRu ] \Big\} \equiv 0,
\end{align}
with the boundary condition $V(t_f,x)=0$. %, which is the Hamilton-Jacobi-Bellman equation of the original optimization problem \eqref{eq:oc1}.
Furthermore, the optimal control $u^*$ is global and unique, since the minimization in \eqref{eq:HJBbilinear} is over a convex (quadratic) function in $u$, and $u^*$ is of the form as expressed in \eqref{eq:ocuk}. %This proves that the convergent solution $u^*$ of the iterative method described in Section \ref{sec:iterative} is a global minimizer to the original Problem \eqref{eq:oc1} given the conditions (i) and (ii). 
\hfill$\Box$

%%%%%%%%%%%%%%%%%%%%%%%%%%%%%%%%%%%%%%%%%%%%%%%%%%%%%%%%%%%%%%%%%%%%%%%%%%
\section{Optimal Control of Bilinear Ensemble Systems}
\label{sec:ensemble}
The iterative method presented in Sections \ref{sec:iterative} and \ref{sec:convergence} can be directly extended to deal with optimal control problems involving a bilinear ensemble system. Consider the minimum-energy control problem for steering a time-invariant bilinear ensemble system, indexed by the parameter $\b$ varying on a compact set $K\subset\mathbb{R}^d$, given by 
\begin{eqnarray}
	\label{eq:bilinear_ensemble}
	\frac{d}{dt}{X(t,\b)}=A(\b)X(t,\b)+B(\b)u+\Big(\sum_{i=1}^m u_i(t) B_i(\b)\Big) X(t,\b),
\end{eqnarray}
where $X=(x_1,\ldots,x_n)^T\in M\subset\mathbb{R}^n$ denotes the state, $\b\in K$, $u:[0,T]\to\mathbb{R}^m$ is the control; the matrices $A(\b)\in\mathbb{R}^{n\times n}$, $B(\b)\in\mathbb{R}^{n\times m}$, and $B_i(\b)\in\mathbb{R}^{n\times n}$, $i=1,\ldots,m$, for $\b\in K$. 
Following the iterative procedure developed in Section \ref{sec:iterative}, we represent the time-invariant bilinear ensemble system in \eqref{eq:bilinear_ensemble} as an iteration equation and formulate the minimum-energy optimal ensemble control problem as
\begin{align}
	\min & \quad J=\frac{1}{2}\int_0^{t_f} (u^{(k)})^T(t)Ru^{(k)}(t)\,dt, \nonumber\\
	\label{eq:oc4}
	{\rm s.t.} & \quad	\frac{d}{dt}{X^{(k)}(t,\b)}=\tilde{A}^{(k-1)}(t,\b)X^{(k)}(t,\b)+\tilde{B}^{(k-1)}(t,\b) u^{(k)}, \tag{P4}\\
	& \quad X^{(k)}(0,\b)=X_0(\b),\quad X^{(k)}(t_f,\b)=X_f(\b), \nonumber
\end{align}
which involves a time-varying linear ensemble system and where we consider the linear ensemble system in a Hilbert space setting; that is, the elements of the matrices $\tilde{A}^{(k-1)}(t,\b)\in\mathbb{R}^{n\times n}$ and $\tilde{B}^{(k-1)}(t,\b)\in\mathbb{R}^{n\times m}$, defined analogously as in \eqref{eq:Ak} and \eqref{eq:Bk}, are real-valued $L_{\infty}$ and $L_2$ functions, respectively, over the space $D=[0,T]\times K$, denoted as $\tilde{A}^{(k-1)}\in L_{\infty}^{n\times n}(D)$ and $\tilde{B}^{(k-1)}\in L_2^{n\times m}(D)$, $X_0,X_f\in L_2^n(K)$, and $R\in\mathbb{R}^{m\times m}\succ 0$.
% which involves a time-varying linear ensemble system and where $X_0,X_f\in L_2^n(K)$ and $R\in\mathbb{R}^{m\times m}\succ 0$; the elements of the matrices $\tilde{A}^{(k-1)}(t,\b)\in\mathbb{R}^{n\times n}$ and $\tilde{B}^{(k-1)}(t,\b)\in\mathbb{R}^{n\times m}$, defined analogously as in \eqref{eq:Ak} and \eqref{eq:Bk}, are real-valued $L_{\infty}^{n\times n}$ and $L_2^{n\times m}$ functions, respectively, and denoted $\tilde{A}^{(k-1)}\in L_{\infty}^{n\times n}$ and $\tilde{B}^{(k-1)}\in L_2^{n\times m}$.

By the variation of constants formula, the ensemble control law that steers the system in \eqref{eq:oc4} between $X_0(\b)$ and $X_f(\b)$ at time $t_f$ satisfies, for each iteration $k$, the integral equation $(L^{(k)}u^{(k)})(\b)=\xi^{(k)}(\b)$,
where
\begin{eqnarray}
	\label{eq:xi}
	\xi^{(k)}(\b)=\Phi^{(k-1)}(0,t_f,\b)X_f(\b)-X_0(\b),
\end{eqnarray}
$\Phi^{(k-1)}(t,0,\b)$ is the transition matrix associated with $\tilde{A}^{(k-1)}(t,\b)$, and where the linear operator $L^{(k)}$ is compact \cite{Li_TAC11} and is defined by
\begin{eqnarray}
	\label{eq:Lk}
	(L^{(k)}u)(\b)=\int_0^T\Phi^{(k-1)}(0,\s,\b)\tilde{B}^{(k-1)}(\b) u(\s)d\s.
\end{eqnarray}

% ============ Theorem 6.1 ==============
\begin{theorem}
	\label{thm:convergence_ensemble}
	Consider the optimal ensemble control problem \eqref{eq:oc4}. Let $(\s_n^{(k)},\mu_n^{(k)},\nu_n^{(k)})$ be a singular system of the operator $L^{(k)}$ defined in \eqref{eq:Lk}. The iterative procedure described according to \eqref{eq:T1} and \eqref{eq:T2} is convergent if the conditions
	\begin{align}
		(i)\ \sum_{n=1}^{\infty}\frac{|\langle\xi^{(k)},\nu_n^{(k)}\rangle|^2}{(\s_n^{(k)})^2}<\infty, \quad (ii)\ \xi^{(k)}\in\overline{\mathcal{R}(L^{(k)})}
	\end{align}
	hold at each iteration $k\in\mathbb{N}$, %$k=1,2,\ldots$, 
	where $\xi^{(k)}$ is defined in \eqref{eq:xi}, $\overline{\mathcal{R}(L^{(k)})}$ denotes the closure of the range space of $L^{(k)}$, and $\langle \xi,\nu\rangle=\int_K \xi^T\nu d\b$ is the inner product defined in $L_2^n(K)$. Furthermore, starting with a feasible initial ensemble trajectory $X^{(0)}(t,\b)$ for Problem \eqref{eq:oc4}, the sequences $X^{(k)}$ and $u^{(k)}$, as in \eqref{eq:x_k} and \eqref{eq:lambda_k}, generated by the iterative method converge to the unique fixed points $X^*$ and $u^*$, respectively.
\end{theorem}

{\it Proof:} Since the conditions (i) and (ii) hold for the operator $L^{(k)}$ at each iteration $k$, the time-varying linear ensemble system in \eqref{eq:oc4} obtained at each iteration $k$ is ensemble controllable, %(or approximate controllable), 
namely, there exists a $u^{(k)}\in L_2^m([0,t_f])$ that steers the ensemble from $X_0(\b)$ to $X_f(\b)$ at time $t_f<\infty$. Moreover, the minimum-energy control that completes this transfer is an infinite weighted sum of the singular functions of $L^{(k)}$, i.e., $u^{(k)}=\sum_{n=1}^{\infty}\frac{1}{\s^{(k)}_n}\langle\xi^{(k)},\nu^{(k)}_n\rangle\mu^{(k)}_n$,
% \begin{eqnarray}
% 	\label{eq:uk_ensemble}
% 	u^{(k)}=\sum_{n=1}^{\infty}\frac{1}{\s^{(k)}_n}\langle\xi^{(k)},\nu^{(k)}_n\rangle\mu^{(k)}_n,
% \end{eqnarray}
and, for any $\e>0$, the truncated optimal control of $u^{(k)}$, i.e.,
\begin{eqnarray}
	\label{eq:uk_N}
	u^{(k)}_N=\sum_{n=1}^{N^{(k)}(\e)}\frac{1}{\s^{(k)}_n}\langle\xi^{(k)},\nu^{(k)}_n\rangle\mu^{(k)}_n,
\end{eqnarray}
drives the ensemble from $X_0(\b)$ to an $\e$-neighborhood of $X_f(\b)$, denoted $\mathcal{B}_\e(X_f)$, satisfying $\|X^{(k)}(t_f,\b)-X_f(\b)\|_2<\e$, where the positive integer $N^{(k)}(\e)$ depends on $\e>0$ \cite{Li_TAC11}. In addition, we denote the optimal trajectories corresponding to the controls $u^{(k)}$ and $u_N^{(k)}$ as $X^{(k)}$ and $X_N^{(k)}$, respectively. Then, according to Theorem \ref{thm:convergence}, the iterative procedure applied to solve for Problem \eqref{eq:oc4} will converge with the convergent optimal control and optimal trajectory pair defined by $(X^*,u^*)$, i.e., $X^{(k)}\rightarrow X^*$ and $u^{(k)}\rightarrow u^*$.

However, the iterations %in the iterative method 
are evolved based on the linear ensemble system formed by the ``truncated trajectory'' $X_N^{(k)}$, given by
\begin{eqnarray}
	\label{eq:truncated_sys}
	\frac{d}{dt}{\hat{X}^{(k+1)}(t,\b)}=\tilde{A}_N^{(k)}(t,\b)\hat{X}^{(k+1)}(t,\b)+\tilde{B}_N^{(k)}(t,\b) \hat{u}^{(k+1)},
\end{eqnarray}
where $\tilde{A}_N^{(k)}(t,\b)=\tilde{A}(X_N^{(k)}(t,\b))$ and $\tilde{B}_N^{(k)}(t,\b)=\tilde{B}(X_N^{(k)}(t,\b))$ depend on $X_N^{(k)}$. %, but not based on the system in Problem \eqref{eq:oc4}. 
Therefore, it requires to show that, at each iteration $k$, the ensemble system as in \eqref{eq:truncated_sys} is ensemble controllable and, furthermore, the optimal control and optimal trajectory obtained based on this iteration equation %\eqref{eq:truncated_sys} 
converge to $X^*$ and $u^*$, respectively, as $k\rightarrow\infty$.

Now, let $L_N^{(k+1)}:L_2^{m}([0,t_f])\rightarrow L_2^n(K)$ be the operator %{\color{red} associated with the truncated control $u_N^{(k-1)}$} 
defined by
\begin{eqnarray}
	\label{eq:LNk}
	\Big(L_N^{(k+1)}u\Big)(\b)=\int_0^T\Phi_N^{(k)}(0,\s,\b)\tilde{B}_N^{(k)}(\s,\b) u(\s)d\s,
\end{eqnarray}
where $\Phi_N^{(k)}(t,0,\b)$ is the transition matrix associated with $\tilde{A}_N^{(k)}(t,\b)$ and $u\in L_2^m([0,t_f])$. Because, by Condition (i),
$u^{(k)}_N$ converges to $u^{(k)}$ uniformly \cite{Li_TAC11} with %such that %\cite{Li_TAC11}
\begin{eqnarray}
	\label{eq:uniform1}
	\|u^{(k)}-u^{(k)}_N\|_2^2 = \sum_{n = N+1}^{\infty } \dfrac{1}{(\s^{(k)}_n)^2}| \langle\xi^{(k)},\nu^{(k)}_n\rangle |^2 \rightarrow 0,\quad\text{as}\quad N\rightarrow\infty,
\end{eqnarray}
and $L^{(k)}$ is compact (see Appendix \ref{appd:linear_ensemble}) so that
\begin{eqnarray}
	\label{eq:uniform2}
	\|L^{(k)}u^{(k)}-L^{(k)}u_N^{(k)}\|_2^2  = \sum_{n = N+1}^{\infty } (\s^{(k)}_n)^2| \langle u^{(k)},u_N^{(k)} \rangle |^2 \rightarrow 0, \quad\text{as}\quad N\rightarrow\infty.
\end{eqnarray}
It follows that $\|X^{(k)}-X^{(k)}_N\|_2^2\rightarrow 0$ as $N\rightarrow\infty$ and, consequently, we have
\begin{eqnarray}
	\label{eq:ANBN}
	\|\tilde{A}^{(k)}-\tilde{A}^{(k)}_N\|_2^2\rightarrow 0,\quad \|\tilde{B}^{(k)}-\tilde{B}^{(k)}_N\|_2^2\rightarrow 0,\quad\text{as}\quad N\rightarrow\infty;
\end{eqnarray}
which implies that at each iteration $k$, the system \eqref{eq:truncated_sys} is also ensemble controllable for sufficiently large $N$.

Next, let $\hat{u}^{(k+1)}$ be the minimum-energy control that steers the system in \eqref{eq:truncated_sys} at each iteration from $X_0(\b)$ to $\mathcal{B}_\e(X_f)$, which is characterized by $(L_N^{(k+1)}\hat{u}^{(k+1)})(\b)=\hat{\xi}^{(k+1)}(\b)$, where $\hat{\xi}^{(k+1)}=\Phi_N^{(k)}(0,t_f,\b)X_f(\b)-X_0(\b)$. Also, we recall that $\Phi_N^{(k)}(t,0,\b)$ and $\Phi^{(k)}(t,0,\b)$ %$\Phi_N^{(k)}(t,0,\b)$ and $\Phi^{(k)}(t,0,\b)$ 
are the respective transition matrices associated with $\tilde{A}_N^{(k)}(\b)\in L_{\infty}^{n\times n}(K)$ and $\tilde{A}^{(k)}(\b)\in L_{\infty}^{n\times n}(K)$, and thus  $\|\Phi^{(k)}(t,0,\b)-\Phi^{(k)}_N(t,0,\b)\|_2^2\rightarrow 0$ as $N\rightarrow\infty$. This guarantees $\|\Phi^{(k)}(0,t,\b)-\Phi^{(k)}_N(0,t,\b)\|_2^2\rightarrow 0$ as $N\rightarrow\infty$ since $\|\Phi_N^{(k)}(t,0,\b)\|$ and $\| \Phi^{(k)}(t,0,\b)\|$ are both bounded. Then, we have
\begin{eqnarray*}
	%\label{eq:udifference}
	\| L^{(k+1)}u^{(k+1)}-L^{(k+1)}_N\hat{u}^{(k+1)}\|_2 = \| \xi^{(k+1)}- \hat{\xi}^{(k+1)} \|_2 \leq  \|\Phi^{(k)}-\Phi_N^{(k)}\|_2 \|X_f\|_2\rightarrow 0
	% \| L^{(k)}u^{(k)}-L^{(k)}_N\hat{u}^{(k)}\|_2 = \| \xi^{(k)}- \hat{\xi}^{(k)} \| \leq \|x_f\| \| \Phi_N^{(k-1)} - \Phi^{(k-1)}\|\rightarrow 0,
\end{eqnarray*}
as $N\rightarrow\infty$ since $X_f\in L_2^n(K)$. This leads to $\|X^{(k+1)}-\hat{X}^{(k+1)}\|_2^2\rightarrow 0$ as $N\rightarrow\infty$, where $\hat{X}^{(k+1)}$ is the trajectory resulting from $\hat{u}^{(k+1)}$. Furthermore, the property $\|\Phi^{(k)}-\Phi^{(k)}_N\|_2\rightarrow 0$, together with \eqref{eq:LNk} and \eqref{eq:ANBN}, gives
\begin{eqnarray}
	\label{eq:LN_L}
	\|L^{(k+1)}_N-L^{(k+1)}\|_2 \rightarrow 0,\quad\text{as}\quad N\rightarrow\infty.
\end{eqnarray}
Because $L^{(k)}(u^{(k)}-\hat{u}^{(k)})=L^{(k)}u^{(k)}-L_N^{(k)}\hat{u}^{(k)}+(L_N^{(k)}-L^{(k)})\hat{u}^{(k)}=\xi^{(k)}-\hat{\xi}^{(k)}+(L_N^{(k)}-L^{(k)})\hat{u}^{(k)}$, 
we obtain, as $N\rightarrow\infty$,
\begin{eqnarray*}
	%\label{eq:udifference}
	\| u^{(k+1)}-\hat{u}^{(k+1)}\|_2 \leq \frac{\|X_f\|_2}{\|L^{(k+1)}\|_2}\|\Phi^{(k)}-\Phi_N^{(k)}\|_2 + \frac{\|\hat{u}^{(k+1)}\|_2}{\|L^{(k+1)}\|_2}\|L^{(k+1)}_N-L^{(k+1)}\|_2\rightarrow 0,
\end{eqnarray*}
by \eqref{eq:LN_L} and by the facts $X_f\in L_2^n(K)$, $\hat{u}^{(k)}\in L_2^m([0,t_f])$, and $\|L^{(k+1)}\|_2 >0$ due to controllability of the system in \eqref{eq:oc4}.

Similar to \eqref{eq:uniform1} and \eqref{eq:uniform2}, we have uniform convergence properties for $\hat{u}_N^{(k)}$, the truncated control of $\hat{u}^{(k)}$, and for $L_N^{(k)}$ such that $\|\hat{u}^{(k)}-\hat{u}^{(k)}_N\|_2^2\rightarrow 0$ and $\| L^{(k)}_N\hat{u}^{(k)}-L^{(k)}_N\hat{u}^{(k)}_N\|_2^2\rightarrow 0$ as $N\rightarrow\infty$. Then, the triangle inequality gives 
\begin{eqnarray}
	\label{eq:uk_u^N}
	\|u^{(k)}-\hat{u}_N^{(k)}\|_2 \leq \|u^{(k)}-\hat{u}^{(k)}\|_2 +\|\hat{u}^{(k)}-\hat{u}_N^{(k)}\|_2 \rightarrow 0,
\end{eqnarray}
and hence
\begin{eqnarray}
	\label{eq:Lu_LNuN}
	% \|L^{(k)}u^{(k)}-L^{(k)}_N\hat{u}^{(k)}_N\|_2^2\rightarrow 0,\quad\text{as}\quad N\rightarrow\infty,
	\|L^{(k)}u^{(k)}-L^{(k)}_N\hat{u}^{(k)}_N\|_2 \leq \|L^{(k)}\|_2\, \|u^{(k)}-\hat{u}_N^{(k)}\|_2 + \|L^{(k)}-L_N^{(k)}\|_2\, \|\hat{u}_N^{(k)}\|_2, %\quad\text{as}\quad N\rightarrow\infty,
\end{eqnarray}
as $N\rightarrow\infty$, which guarantees that, at each iteration $k$, the trajectory $\hat{X}_N^{(k)}$ converges to $X^{(k)}$, which result from $\hat{u}^{(k)}_N$ and $u^{(k)}$, respectively, i.e.,
\begin{eqnarray}
	\label{eq:x_N}
	\|X^{(k)}-\hat{X}_{N}^{(k)}\|_2^2\rightarrow 0, \quad \text{as}\quad N\rightarrow \infty.
	% \| \hat{x}_{N}^{(k)}-x^{(k)}\|_2^2\rightarrow 0, \text{ as } N\rightarrow \infty,
\end{eqnarray}
In addition, since $X^{(k)}\rightarrow X^*$ and $u^{(k)}\rightarrow u^*$ as $k\rightarrow\infty$, we have 
\begin{eqnarray}
	\label{eq:uu*}
	\|\hat{u}_N^{(k)}-u^*\|_2 \leq \|\hat{u}_N^{(k)}-u^{(k)}\|_2 +\|u^{(k)}-u^*\|_2 \rightarrow 0
\end{eqnarray}
as $k,~N\rightarrow\infty$ by \eqref{eq:uk_u^N}, as well as $\| L^{(k)}-L^*\|_2\rightarrow 0$ as $k\rightarrow\infty$, where $L^*$ is the operator defined with respect to the convergent solutions $X^*$ and $\l^*$, given by
\begin{eqnarray}
	\label{eq:L}
	(L^*u)(\b)=\int_0^T\Phi^*(0,\s,\b)\tilde{B}\big(X^*(t,\b)\big)u(\s)d\s,
\end{eqnarray}
in which $\Phi^*(t,0,\b)$ is the transition matrix associated with $\tilde{A}(X^*(t,\b))$. This then gives 
\begin{eqnarray}
	\label{eq:LuLu*}
	\|L^{(k)}u^{(k)}-L^*u^*\|_2 \leq \|L^{(k)}\|_2 \|u^{(k)}-u^*\|_2 +\|L^{(k)}-L^*\|_2 \|u^*\|_2 \rightarrow 0,
\end{eqnarray}
as $k\rightarrow\infty$. Finally, combining \eqref{eq:Lu_LNuN} and \eqref{eq:LuLu*} and applying the triangle inequality yield 
$$\|L^*u^*-L^{(k)}_N\hat{u}^{(k)}_N\|_2\rightarrow 0 \text{\quad as \quad} k,N\rightarrow\infty ,$$ 
which implies that $\|X^*-\hat{X}_N^{(k)}\|\rightarrow 0$ as $k,N\rightarrow\infty$. This together with \eqref{eq:uu*} concludes the convergence of the sequences $\{\hat{u}_N^{(k)}\}$ and $\{\hat{X}_N^{(k)}\}$, generated by the iterative method, to $u^*$ and $X^*$, respectively, i.e., the minimum-energy ensemble control law and the optimal ensemble trajectory that satisfy the necessary optimality condition. \hfill$\Box$

%%%%%%%%%%%%%%%%%%%%%%%%%%%%%%%%%%%%%%%%%%%%%%%%%%%%%
\section{Examples and Numerical Simulations}
\label{sec:examples}
In this section, we apply the developed iterative algorithm to solve optimal control problems involving single and ensemble bilinear systems, including the well-known Bloch system that models the evolution of two-level quantum systems \cite{Li_PNAS11}. 
Ensemble control of Bloch systems is a key to many applications in quantum control, such as nuclear magnetic resonance %(NMR) 
spectroscopy and imaging (MRI), quantum computation and quantum information processing \cite{Li_PRA06, Mabuchi05}, as well as quantum optics \cite{Silver85}.

\begin{figure}[t]
 	\centering
 	\subfigure[The convergent optimal control]{\includegraphics[width=0.49\columnwidth]{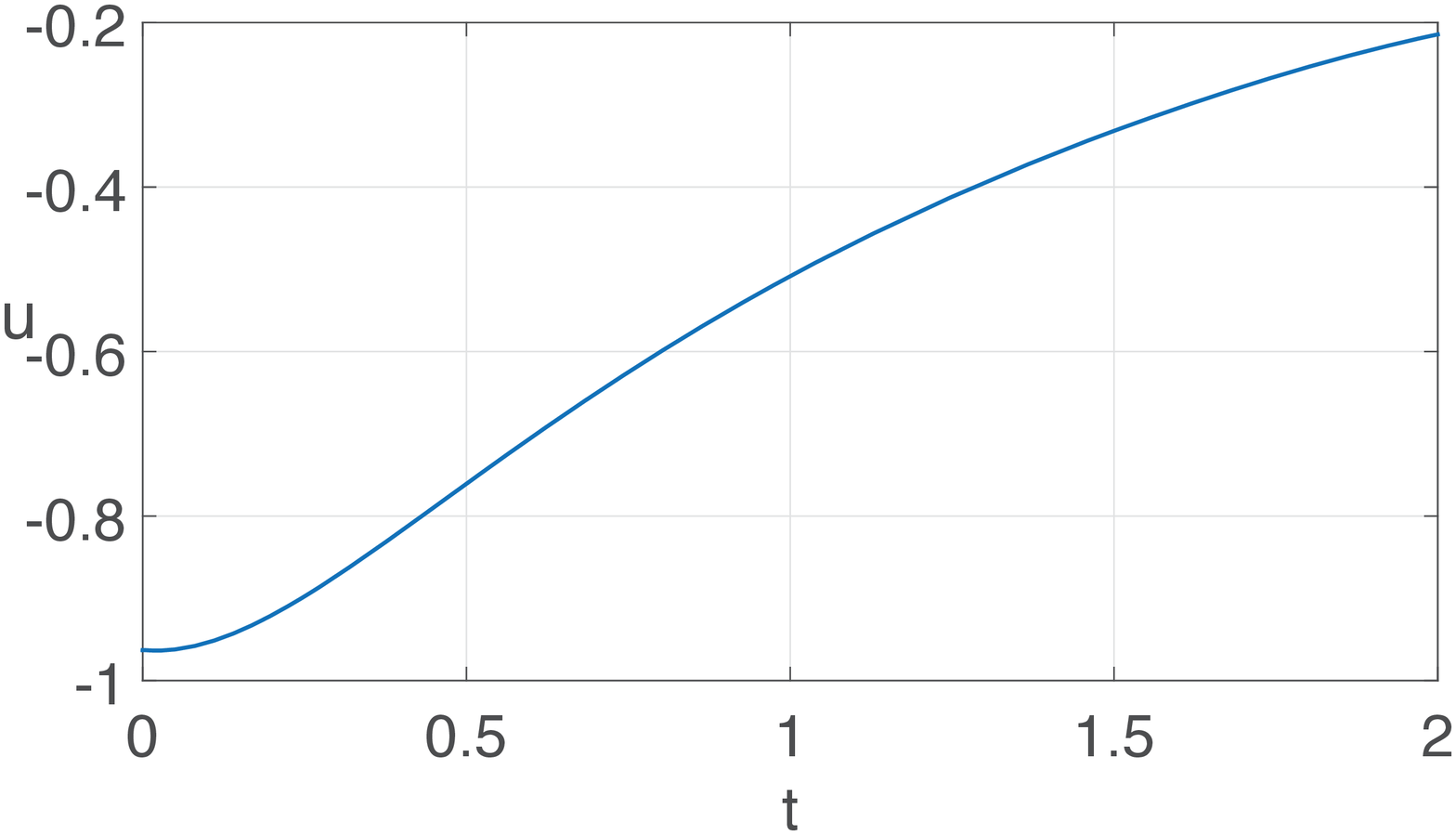}
	\label{fig:PControl}}
 	\subfigure[The optimal trajectory]{\includegraphics[width=0.49\columnwidth]{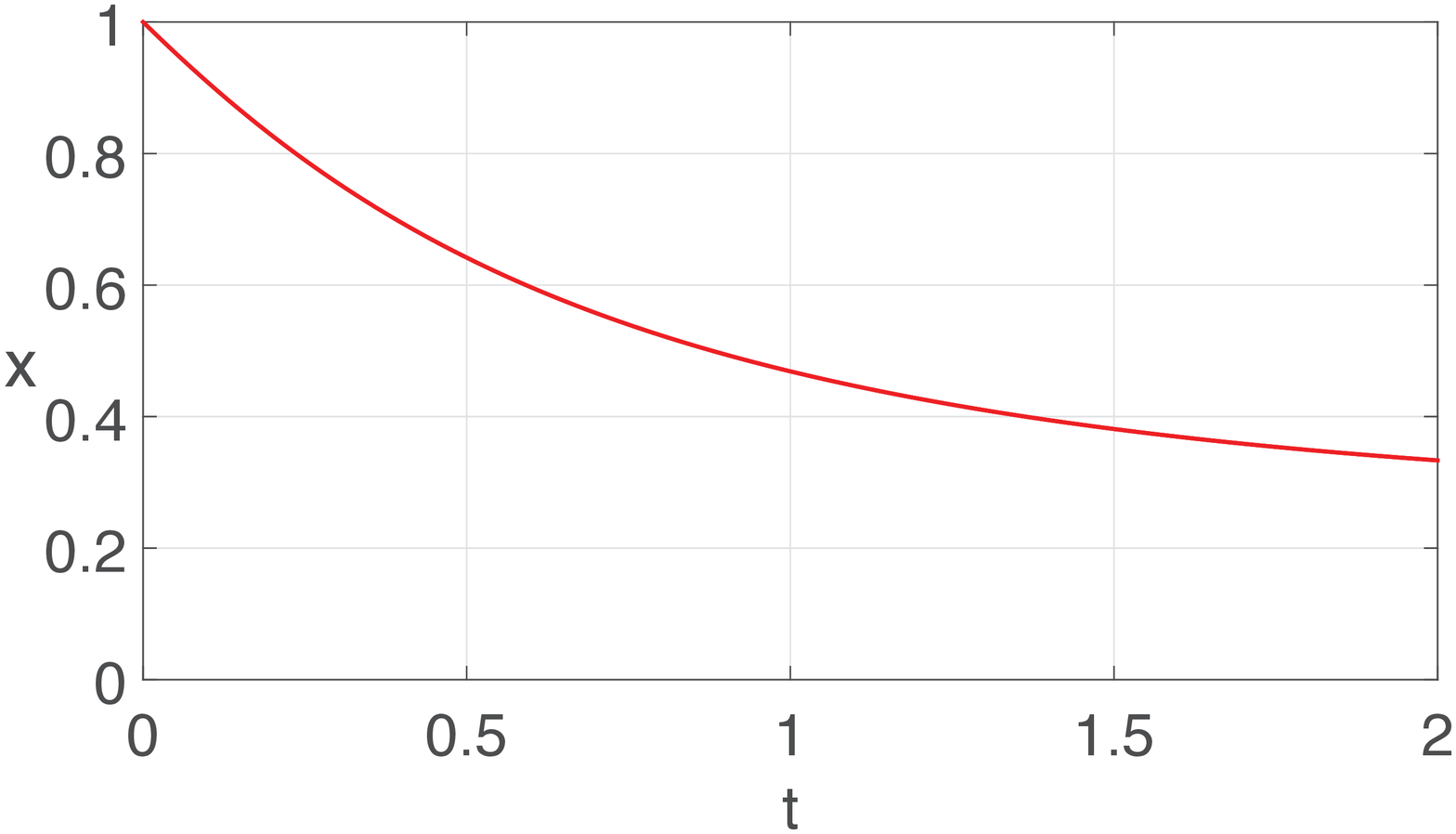}
	\label{fig:PState}} 
 	\caption{\subref{fig:PControl} The optimal control that steers the bilinear system in Example \ref{ex:population} %\eqref{eq:psys} 
from $x(0)=1$ to $x(2)=1/3$, while minimizing %the objective functional 
$J=\int_0^2 (x^2+u^2)dt$. 
	\subref{fig:PState} The optimal trajectory following the optimal control shown in \subref{fig:PControl}.} %\vspace{-24pt}
	\label{fig:PopulationSingle}
 \end{figure}

% =============== Example 1 ====================
\begin{example}[Population Control in Socioeconomics] 
\label{ex:population}
\rm We consider a simple but representative bilinear system arising from the field of socioeconomics, which models the dynamics of population growth, simplified based on the Gibson's population transfer model \cite{Gibson1972}, given by $\frac{dx}{dt} = ux$,
% \begin{align}
% 	\label{eq:psys}
% 	\frac{dx}{dt} = ux,
% \end{align}
where $x\in\mathbb{R}$ represents the number of domestic laborers and $u\in\mathbb{R}$ denotes the attractiveness for immigration multiplier. Putting this into the canonical form as presented in \eqref{eq:linear}, we have $A=0$, $B=0$, and $N=1$. We consider the design of the optimal control $u$ for reducing two-third of the domestic laborer population, i.e., from $x(0)=1$ to $x(t_f)=1/3$, in $t_f=2$, which minimizes the cost functional $J=\int_0^2 (x^2+u^2)dt$. The optimal control obtained by the iterative method is shown in Figure \ref{fig:PControl} and the resulting optimal trajectory is displayed in Figure \ref{fig:PState}.
\end{example}

% ===================== Figure 1 =======================

% ===================== Figure 2 =======================
\begin{figure}[t]
 	\centering
 	\subfigure[The convergent minimum-energy control and the resulting optimal trajectory]{\includegraphics[width=0.55\columnwidth]{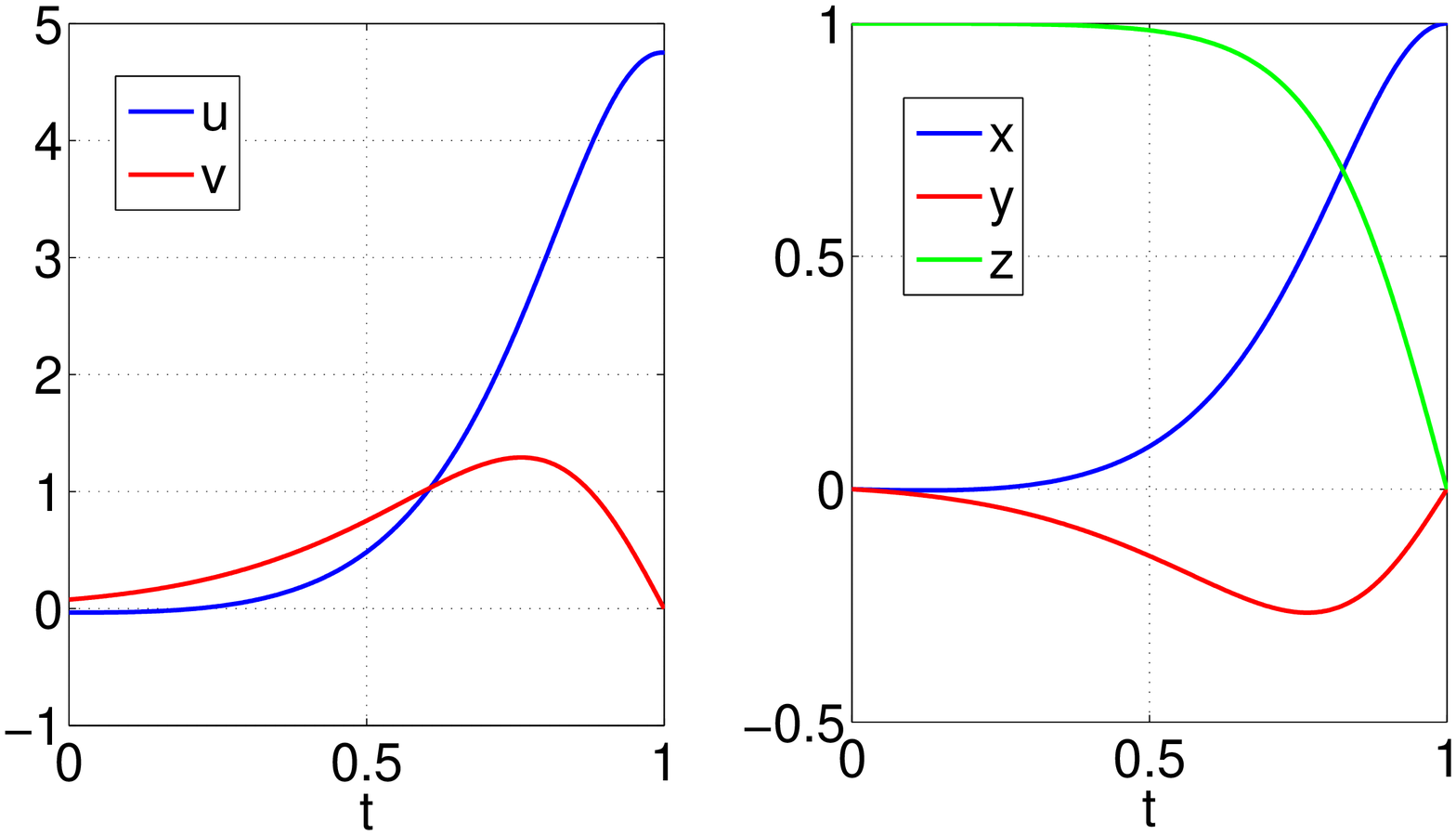} \label{fig:singlecontrol}}
 	\subfigure[The optimal trajectory on the Bloch sphere]{\includegraphics[width=0.38\columnwidth]{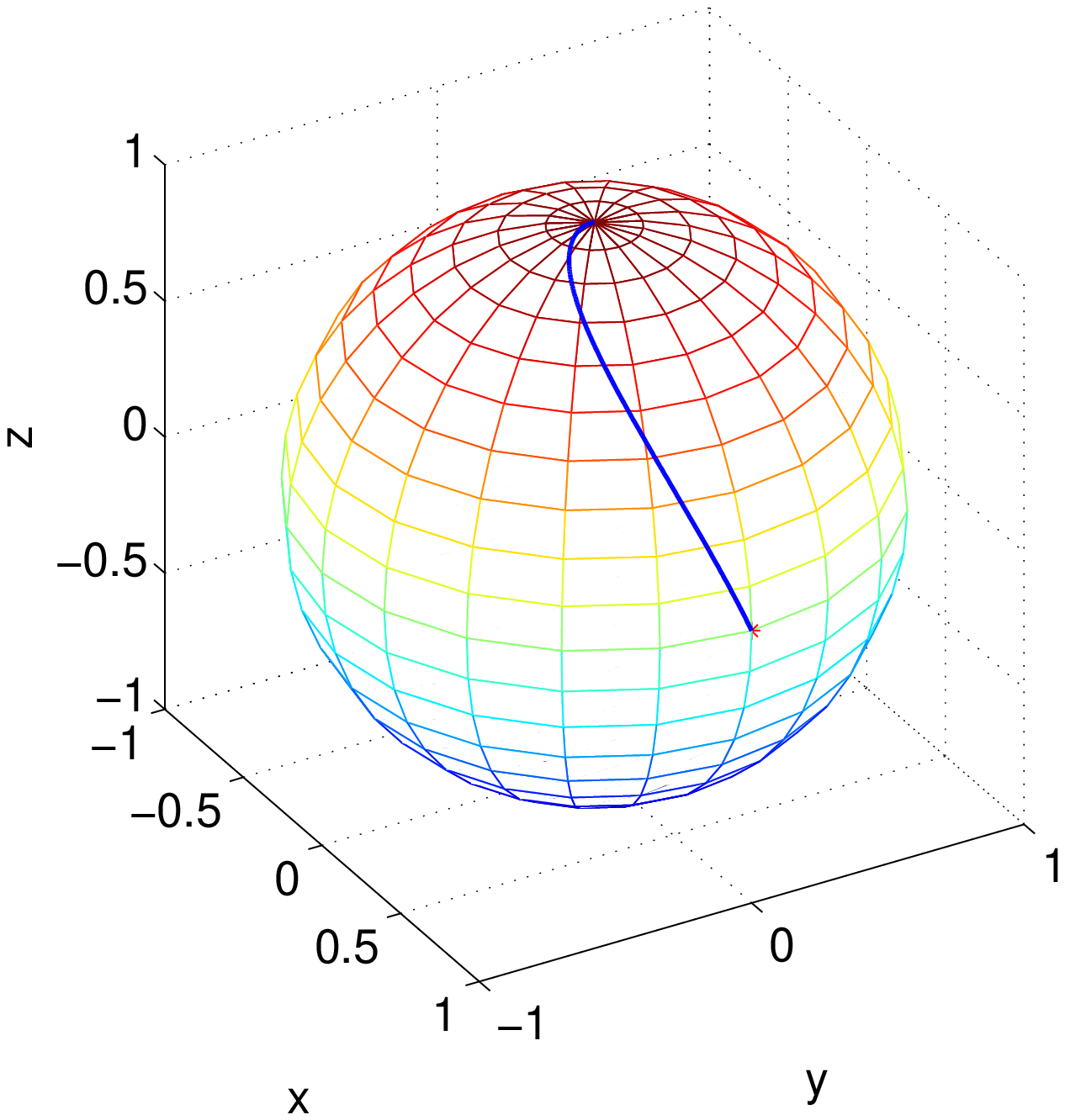}\label{fig:singletrajectory}} 
 	\caption{\subref{fig:singlecontrol} The minimum-energy excitation ($\pi/2$) pulse (left panel) that steers the Bloch system in \eqref{eq:blochsys} with $\omega=0.5$ from $x_0 = (0,0,1)^T$ to $x_f=(1,0,0)^T$ at $t_f=1$ and the resulting optimal trajectory (right panel). \subref{fig:singletrajectory} The optimal trajectory on the Bloch sphere}. %\vspace{-24pt}
 \label{fig:blochsingle}
 \end{figure}

% =============== Example 2 ====================
% \subsection{Steering a Single Bloch System}
\begin{example}[Excitation of a Two-Level System]
\label{ex:bloch}
\rm A canonical example of optimal control of bilinear ensemble systems in quantum control is the optimal pulse design for the excitation of a collection of two-level systems \cite{Li_PNAS11}, in which the dynamics of a quantum %spin 
ensemble obeys the Bloch equations, and optimal pulses (controls) that steer the ensemble between states of interest are pursued. The Bloch equations form a bilinear control system evolving on the special Lie group SO(3), given by 
\begin{align}
	\label{eq:blochsys}
	\frac{d}{dt} \begin{bmatrix}
	x_1 \\ x_2 \\ x_3
	\end{bmatrix} = \begin{bmatrix}
	0 & -\omega & u_1 \\
	\omega & 0 & -u_2 \\
	-u_1 & u_2 & 0 \end{bmatrix}
	\begin{bmatrix}
		x_1 \\ x_2 \\ x_3
	\end{bmatrix},
\end{align}
where $x=(x_1,x_2,x_3)^T$ denotes the bulk magnetization of the spins, $\w$ denotes the Larmor frequency of the spins, and $u_1$ and $u_2$ are the radio-frequency fields applied on the $y$ and the $x$ direction, respectively \cite{Li_TAC09}. A common control task is to drive the system from the equilibrium state $x_0 = (0,0,1)^T$ to an excited state on the transverse plane, e.g., $x_f=(1,0,0)^T$, and, in particular, achieving the desired state transfer with minimum-energy is of practical importance \cite{Li_PNAS11}.

Here, we consider exciting a spin system with the Larmor frequency $\w=0.5$, and first rewrite the Bloch system in the canonical form as presented in \eqref{eq:linear} with 
$$A = \begin{bmatrix}
	0 & -\omega & 0 \\
	\omega & 0 & 0 \\
	0 & 0 & 0
	\end{bmatrix}, \quad B=0;
	 \quad N_1 = \begin{bmatrix}
	0 & 0 \\
	0 & 0 \\
	-1 & 0
	\end{bmatrix}, \quad
	N_2 = \begin{bmatrix}
	0 & 0 \\
	0 & 0 \\
	0 & 1
	\end{bmatrix}, \quad
	N_3 = \begin{bmatrix}
	1 & 0 \\
	0 & -1 \\
	0 & 0
	\end{bmatrix},$$
and apply the iterative method described in Section \ref{sec:iterative} to find the minimum-energy control. Figure \ref{fig:singlecontrol} illustrates the convergent minimum-energy control that steers the spin system from $x_0$ to $x_f$ at $t_f=1$ and minimizes $J=\int_0^1 (u_1^2+u_2^2) dt$, and the resulting trajectory is shown in Figure \ref{fig:singletrajectory}. This optimal control and the trajectory converge in 17 iterations, starting with an initial trajectory $x^{(0)}$ with endpoints $x_0$ and $x_f$ with the least distance, given the stopping criterion $\|x(t_f)-x_f\|\leq 10^{-5}$.
\end{example}

\begin{figure}[t]
	\centering
 	\subfigure[The minimum-energy ensemble control law]{\includegraphics[width=0.45\columnwidth]{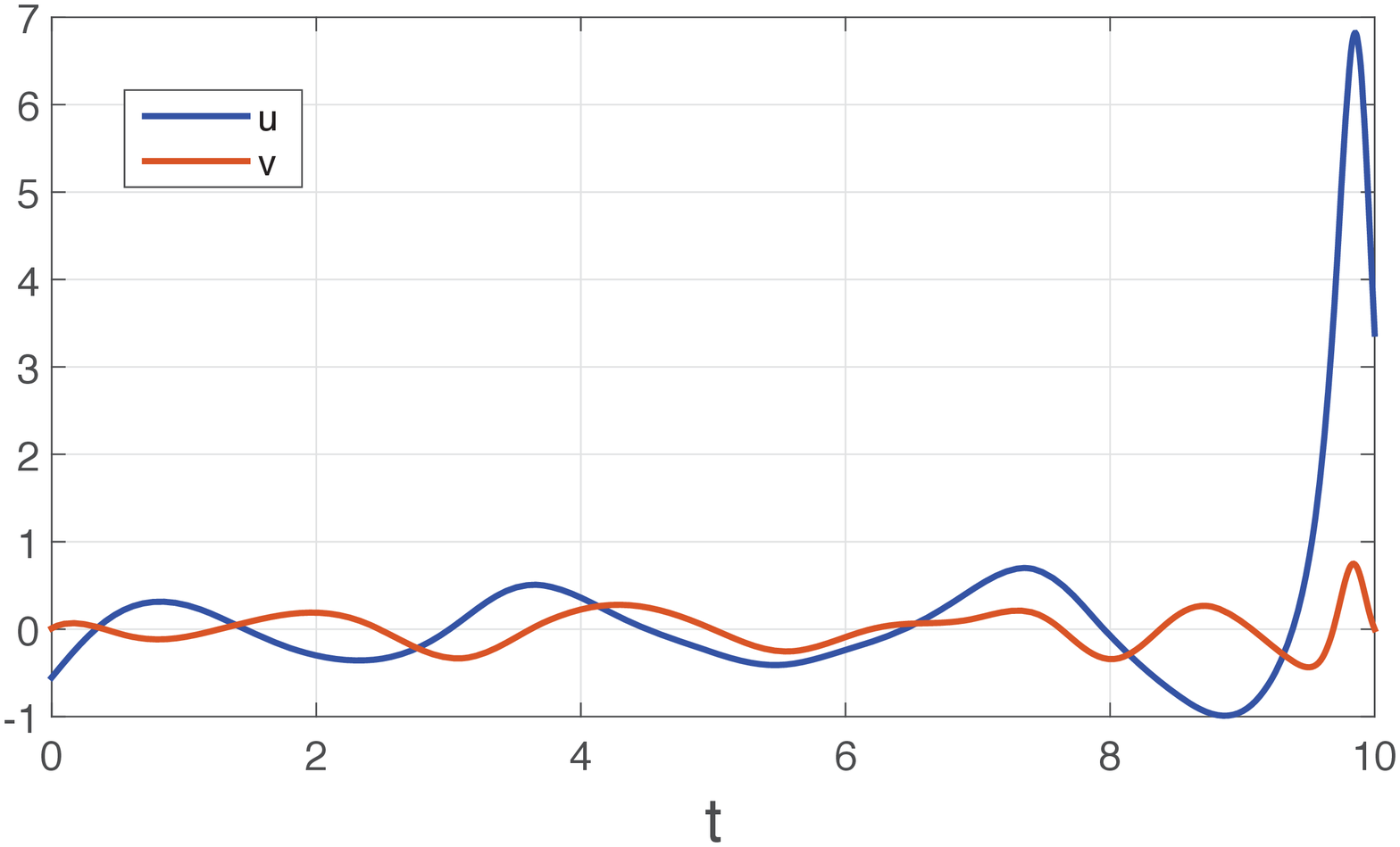} \label{fig:ensemblecontrol}}
 	\subfigure[The $x$-components of the final states of the ensemble]{\includegraphics[width=0.48\columnwidth]{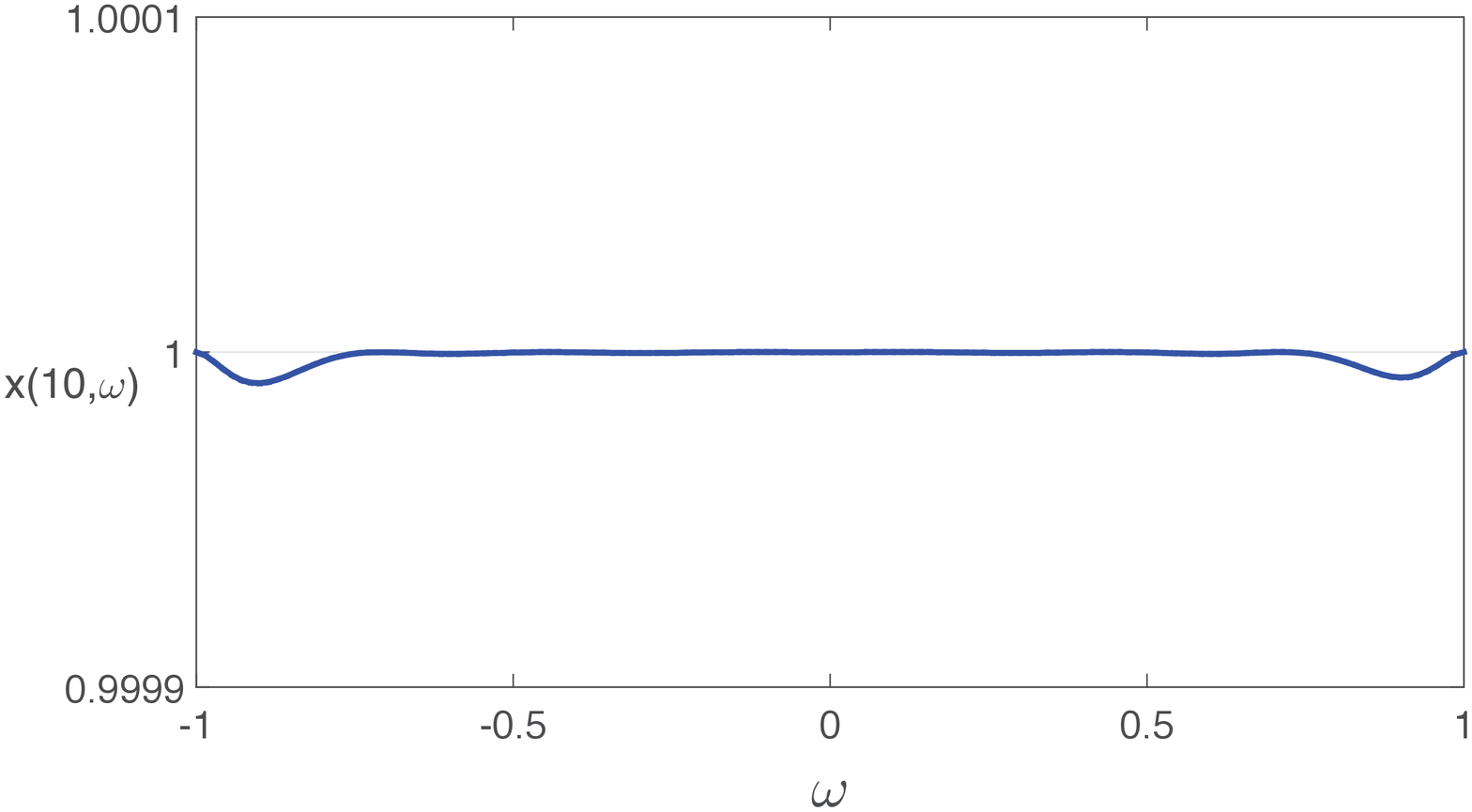}\label{fig:ensembletrajectory}} 
 	\caption{\subref{fig:ensemblecontrol} The minimum-energy ensemble control that steers an ensemble of Bloch systems with $\w\in[-1,1]$ from $X_0(\w) = (0,0,1)^T$ to a neighborhood of $X_f(\w)=(1,0,0)^T$. The weighted matrix $R=I$. \subref{fig:ensembletrajectory} The $x$-components of the final states $X(10,\w)$ for 141 spin systems with their frequencies uniformly spaced within $[-1,1]$ following the minimum-energy control displayed in \subref{fig:ensemblecontrol}.}
 \label{fig:blochensemble}
 \end{figure}

% =================== Example 3 =======================
\begin{example}[Excitation of an Ensemble of Two-Level Systems]
\label{ex:bloch_ensemble}
\rm Here, we apply the iterative method to design a minimum-energy broadband excitation ($\pi/2$) pulse that steers an ensemble of spin systems modeled in \eqref{eq:blochsys} with the ensemble state defined as $X(t,\w)=(x(t,\w),y(t,\w),z(t,\w))^T$ for $\w\in [-1,1]$ from $X(0,\w)=(0,0,1)^T$ to $X(t_f,\w)=(1,0,0)^T$, where $t_f=10$ is the pulse duration. At each iteration $k$, the minimum-energy ensemble control law is calculated using an singular-value-decomposition (SVD) based algorithm \cite{Li_ACC12_SVD}, by which the input-to-state operator $L^{(k)}_N$ as in \eqref{eq:LNk} that characterizes the evolution of the system dynamics is approximated by a matrix of finite rank, due to the compactness of this operator. Then, the singular values, $\s_n^{(k)}$, and the singular vectors, $\mu_n^{(k)}$ and $\nu_n^{(k)}$, are calculated using SVD to synthesize the optimal ensemble control expressed in \eqref{eq:uk_N} \cite{Li_ACC12_SVD}. The convergent minimum-energy ensemble control law, i.e., the minimum-energy broadband $\pi/2$ pulse, is illustrated in Figure \ref{fig:ensemblecontrol}, and the performance, i.e., the $x$-component of the final state $X(t_f,\w)$, $t_f=10$, for 141 spin systems is shown in Figure \ref{fig:ensembletrajectory}. The iterative algorithm converges in 152 iterations % to get $u^*=u^{(152)}$, 
given the stopping criterion $\|X(t_f,\w)-X_f(\w)\|\leq 10^{-5}$. %{\color{red} (It would sound better to also show the computational time on your PC, for all of the examples)}.

\end{example}

%================= Figure 3 ===========================   

%%%%%%%%%%%%%%%%%%%%%%%%%%%%%%%%%%%%%%%%%%%%%%%%%%%%%%%%%%%%%%%%%%%%%%%%%%%%%%%%
\section{Conclusion} %S AND FUTURE WORKS}
We develop an iterative method for solving fixed-endpoint optimal control problems involving time-invariant bilinear and bilinear ensemble systems. We analyze the convergence of the iterative procedure by using the contraction mapping and the fixed-point theorem. The central idea of our approach is to represent the time-invariant bilinear ensemble system as a time-varying linear ensemble system and then to show, in an iterative manner, that the optimal control of the original bilinear (ensemble) system is the convergent optimal control of the associated linear (ensemble) system. In addition, we illustrate the condition for global optimality of the convergent solution. Finally, we demonstrate the effectiveness and applicability of the constructed iterative method using several examples involving the control of population growth in socioeconomics and the design of broadband pulses for exciting an ensemble of two-level systems, which is essential to many applications in quantum control.

%%%%%%%%%%%%%%%%%%%%%%%%%%%%%%%%%%%%%%%%%%%%%%%%%%%%%%%%%%%%%%%%%%%%%%%%%%
\section{Appendix}
\subsection{Sweep method and the notion of flow mapping}%for fixed-endpoint optimal control problems}
We illustrate the idea of the developed Sweep method for dealing with fixed-endpoint optimal control problems presented in Section \ref{sec:sweep} using the notion of flow mapping \cite{Schaettler13} and show the connection between the non-singularity of the matrix $P$ defined in \eqref{eq:pode} and the controllability of the system in \eqref{eq:oc2}.

% Here, we explain the iterative procedure using the flowing mapping \cite{Schaettler13} and derive the related ordinary differential equations (ODE's). Then, we show the relationship between the non-singularity of the matrix $P$ and the controllability of the system \eqref{eq:oc2} and claim these two are actually equivalent.

% ---------------------------------------------------
\subsubsection{Flow mapping in optimal control} %Sweep method for fix-point problems}
\label{appd:sweep}
% We introduce the Sweep method using the concepts of flow mapping. 
Consider the optimal control problem parameterized by $p$ given by,
\begin{align}
	\label{eq:oc_flow}
	\min &\quad J=\int_0^{t_f} L(t,x(t,p),u(t,p)) \, dt, \nonumber\\
	{\rm s.t.} &\quad \dot{x}(t,p)=f(t,x(t,p),u(t,p)), \tag{P5}
\end{align}
where $p\in\mathbb{R}$ represents the perturbation of extremals. Let $\mathcal{E}$ be a $C^r$-parameterized family of extremals for Problem \eqref{eq:oc_flow}, and suppose that $F:E\rightarrow G$ is the flow restricted on a subset $E$ of the $(t,p)$-space, which is a $C^{1,r}$-diffeomorphism onto an open subset $G\subset \mathbb{R}\times \mathbb{R}^n$ of the $(t,x)$-space. Let $C$ be the cost-to-go function for Problem \eqref{eq:oc_flow}, then the value function $V^{\mathcal{E}}:G\rightarrow \mathbb{R}$ of $\mathcal{E}$ defined by
$V^{\mathcal{E}} = C\circ F^{-1}$,
is continuously differentiable in $(t,x)$ and $r$-times continuously differentiable in $x$ for any fixed $t$ (see Theorem 5.2.1 in \cite{Schaettler13}). In addition, the function $u_*:G\rightarrow \mathbb{R}$ defined by $u_* = u\circ F^{-1}$
% \begin{align*}
% 	u_* = u\circ F^{-1}, 
% \end{align*}
is an admissible feedback control that is continuous and $r$-times continuously differentiable in $x$ for any fixed $t$. The diagrams below illustrate the relation of the mappings defined above.

% ============ Diagram ===========
\begin{center}
\qquad \xymatrix{
E \ar[d]^F \ar[r]^C &{\mathbb{R}}\\
G \ar[ru]_{V^{\mathcal{E}}}}
\qquad \qquad \qquad \xymatrix{
E \ar[d]^F \ar[r]^u &{\mathbb{R}}\\
G \ar[ru]_{u_*}}
\end{center}

\noindent Together, the pair $(V^{\mathcal{E}},u_*)$ is a classical solution to the Hamilton-Jacobi-Bellman equation, and the following identities hold for all $p$,
\begin{align*}
\dfrac{\partial V^{\mathcal{E}}}{\partial t}(t, x(t,p)) & = -H(t,\l(t,p),x(t,p),u(t,p)), \\
\dfrac{\partial V^{\mathcal{E}}}{\partial x}(t, x(t,p)) & = \l(t,p),
\end{align*}
where $H$ is the Hamiltonian associated with Problem \eqref{eq:oc_flow} and $\lambda$ is the co-state of $x$.
$V^{\mathcal{E}}$ is $(r+1)$-times continuously differentiable in $x$ on $G$ because $\mathcal{E}$ is nicely $C^r$-parameterized, and then we have $\dfrac{\partial^2 V^{\mathcal{E}}}{\partial x^2}(t, x(t,p)) = \dfrac{\partial \l}{\partial p}(t,p)\Big(\dfrac{\partial x}{\partial p}(t,p)\Big)^{-1}$,
% \begin{align*}
% 	\dfrac{\partial^2 V^{\mathcal{E}}}{\partial x^2}(t, x(t,p)) = \dfrac{\partial \l}{\partial p}(t,p)\Big(\dfrac{\partial x}{\partial p}(t,p)\Big)^{-1},
% \end{align*}
provided that $\dfrac{\partial x}{\partial p}(t,p)$ is invertible.

% ---------------------------------------
\subsubsection{Sweep method derived based on the flow mapping}
Now, consider the parameterized optimal control problem associated with Problem \eqref{eq:oc2} %associated with the system in in Problem \eqref{eq:oc2}, 
given by 
\begin{align}
	\min &\quad J=\frac{1}{2}\int_0^{t_f} \Big[x^T(t,p)Qx(t,p)+u^T(t,p)Ru(t,p)\Big] \, dt, \nonumber\\
	\label{eq:oc6}
	{\rm s.t.} &\quad 
	\dot{x}(t,p)=\tilde{A}x(t,p)+\tilde{B}u(t,p) \tag{${\rm P2}_p$}\\  
	&\quad x(0,p)=x_0, \quad x(t_f,p)=x_f, \nonumber
\end{align}
we then have, for this parameterized LQR problem,
\begin{align*}
	\dfrac{d}{dt}\Big(\dfrac{\partial x}{\partial p}(t,p)\Big) &= \tilde{A}\dfrac{\partial x}{\partial p}(t,p)-\tilde{B}R^{-1}\tilde{B}^T\dfrac{\partial \l}{\partial p}(t,p), \\ 
	\dfrac{d}{dt}\Big(\dfrac{\partial \l}{\partial p}(t,p)\Big) &= -Q\dfrac{\partial x}{\partial p}(t,p)-\tilde{A}^T\dfrac{\partial \l}{\partial p}(t,p).
\end{align*}
Defining
\begin{eqnarray}
	\label{eq:Delta}
	\Delta (t,p) = \dfrac{\partial \l}{\partial p}(t,p) - K(t,p)\dfrac{\partial x}{\partial p}(t,p)-S(t,p)\dfrac{\partial \nu}{\partial p}(p)
\end{eqnarray}
with $\Delta (t_f,p) = 0$, where $K$ and $S$ satisfy the matrix differential equations \eqref{eq:kRiccati} and \eqref{eq:s} with the terminal conditions $K(t_f,p)=0\in \mathbb{R}^{n\times n}$ and $S(t_f,p)=I\in \mathbb{R}^{n\times n}$, respectively; \eqref{eq:Delta} yields, %with the above,
\begin{align*}
\dot{\Delta}(t,p) & = \dfrac{d}{dt}\Big(\dfrac{\partial \l}{\partial p}\Big) - \dot{K}\dfrac{\partial x}{\partial p} - K\dfrac{d}{dt}\Big(\dfrac{\partial x}{\partial p}\Big) - \dot{S}\dfrac{\partial \nu}{\partial p}(p) \\
& = -Q\dfrac{\partial x}{\partial p} - \tilde{A}^T\dfrac{\partial \l}{\partial p} - (-Q-K\tilde{A}-\tilde{A}^TK+K\tilde{B}R^{-1}\tilde{B}^TK)\dfrac{\partial x}{\partial p} \\
& - K \Big(\tilde{A}\dfrac{\partial x}{\partial p}-\tilde{B}R^{-1}\tilde{B}^T\dfrac{\partial \l}{\partial p}\Big) + [\tilde{A}^T-K^T\tilde{B}R^{-1}\tilde{B}^T]S\dfrac{\partial \nu}{\partial p}(p) \\
& = -[\tilde{A}^T-K^T\tilde{B}R^{-1}\tilde{B}^T]\Delta(t,p).
\end{align*}
This gives $\Delta(t,p) \equiv 0$, since $\Delta (t_f,p) = 0$, and hence guarantees, from \eqref{eq:Delta}, that $\dfrac{\partial \l}{\partial p}(t,p) = K(t,p)\dfrac{\partial x}{\partial p}(t,p)+S(t,p)\dfrac{\partial \nu}{\partial p}(p)$,
% \begin{align*}
% 	\dfrac{\partial \l}{\partial p}(t,p) = K(t,p)\dfrac{\partial x}{\partial p}(t,p)+S(t,p)\dfrac{\partial \nu}{\partial p}(p),
% \end{align*}
which we adopted to define \eqref{eq:lam}.

In order to fulfill the terminal condition $x(t_f,p)=x_f$ at time $t_f$, we introduce the auxiliary variables $O(t,p),~P(t,p)\in\mathbb{R}^{n\times n}$， and set
 \begin{eqnarray}
	\label{eq:end}
 	x_f = O(t,p)x(t,p)+P(t,p)\nu(p).
 \end{eqnarray}
Clearly, at $t=t_f$, we need $O(t_f,p)=I$ and $P(t_f,p)=0$. Taking the time derivative on both sides of \eqref{eq:end}, we get $0=\dot{O}(t,p)x(t,p)+O(t,p)\dot{x}(t,p)+\dot{P}(t,p)\nu(p)$, which results in
 % \begin{align*}
 %  \dot{O}(t,p)x(t,p)+O(t,p)[\tilde{A}x(t,p)-\tilde{B}R^{-1}\tilde{B}^T(K(t,p)x(t,p)+S(t,p)\nu(p))]+\dot{P}(t,p)\nu(p)=0.
 %  \end{align*}
 % We rewrite it in the form, 
\begin{align*}
	\big\{\dot{O}(t,p)+O(t,p)[\tilde{A}-\tilde{B}R^{-1}\tilde{B}^TK(t,p)]\big\}x(t,p)+[\dot{P}(t,p)-O(t,p)\tilde{B}R^{-1}\tilde{B}^TS(t,p)]\nu(p)=0.
\end{align*}
Since it holds for all $x$ and $\nu$, we have
\begin{align}
	\label{eq:Odot}
	\dot{O}(t,p)+O(t,p)[\tilde{A}-\tilde{B}R^{-1}\tilde{B}^TK(t,p)]=0 \quad &\text{ with}\quad O(t_f,p)=I, \\
	\label{eq:Pdot}
	\dot{P}(t,p)-O(t,p)\tilde{B}R^{-1}\tilde{B}^TS(t,p)=0 \quad &\text{ with}\quad P(t_f,p)=0.
\end{align}
Observe that, from \eqref{eq:Odot}, $\dot{O}^T(t,p)$ satisfies the same equation, for all $p$, as $S(t)$ in \eqref{eq:s} with the terminal condition $S(t_f,p)=I$; hence, $O(t,p)=S^T(t,p)$ for $t\in [0,t_f]$. Then, we can rewrite \eqref{eq:Pdot} as
\begin{align}
	\dot{P}(t,p)-S^T(t,p)\tilde{B}R^{-1}\tilde{B}^TS(t,p)=0
\end{align}
for $t\in [0,t_f]$ with the terminal condition $P(t_f,p)=0$.

The multiplier associated with the terminal constraint can be expressed, by \eqref{eq:end}, as $\nu(p)=\big[P(t,p)\big]^{-1} \big[x_f-S^T(t,p)x(t,p)\big]$, provided $P(t,p)$ is invertible for $t\in[0,t_f]$. Note this condition holds for all $t\in [0,t_f]$, i.e., $\nu(p)$ is a constant for a fixed $p$; hence plugging in $t=0$, we have $\nu(p)$ in terms of the initial and terminal states, given by
\begin{equation}
	\label{eq:nu_p}
	\nu(p)=\big[P(0,p)\big]^{-1} \big[x_f-S^T(0,p)x_0\big].
\end{equation}

% ------------------------------------
\subsubsection{Controllability of the system in \eqref{eq:oc6} and non-singularity of the matrix $P$}
\label{appd:Pnonsingular}
The following two lemmas will illustrate that controllability of the system in \eqref{eq:oc6} guarantees non-singularity of the matrix $P$, so that $\nu(p)$ in \eqref{eq:nu_p} is well defined.

% Before proving the equivalence of the two concepts, we introduce two lemmas which could help us illustrate the idea clearly.

% ==========================================
\begin{lemma}
	\label{lem:Psingular}
	The matrix $P(\tau,p)$ is singular if and only if there exists a nontrivial solution $\mu =\mu(t,p)$ of the linear adjoint equation
	$\dot{\mu} = -\mu \dfrac{\partial f}{\partial x} (t,x(t,p),u(t,p))$
	with the terminal condition $\mu(t_f,p)$ that is perpendicular to the terminal manifold at $x(t_f,p)$ such that 
	$\mu (t,p)\dfrac{\partial f}{\partial u} (t,x(t,p),u(t,p)) \equiv 0$
	on the interval $[\tau, t_f]$, where $\dot{x}(t,p) = f(t,x(t,p),u(t,p))$ \cite{Schaettler13}.
\end{lemma}

\begin{lemma}
	\label{lem:Pcontrollable}
	A time-varying linear system $\dot{x}=A(t)x+B(t)u$ is controllable over an interval $[\tau, t_f]$ if and only if for every nontrivial solution $\mu$ of the adjoint equation $\dot{\mu} = -\mu A(t)$, the function $\mu(t)B(t)$ does not vanish identically on the interval $[\tau, t_f]$. It is completely controllable if this holds for any subinterval $[\tau, t_f]$ \cite{Schaettler13}.
\end{lemma}

By Lemma \ref{lem:Pcontrollable}, the linear system $\dot{x}=\tilde{A}(t)x+\tilde{B}(t)u$ in Problem \eqref{eq:oc6} is controllable over an interval $[t, t_f]$ if and only if for every nontrivial solution $\mu$ of the adjoint equation $\dot{\mu} = -\mu \tilde{A}(t)$, the function $\mu(t)\tilde{B}(t)$ does not vanish identically on the interval $[t, t_f]$, which by Lemma \ref{lem:Psingular} is equivalent to the non-singularity of $P(t,p)$, $\forall t \in [0, t_f]$. In particular, for $t=0$, we have that $P(0,t_f)$ is invertible if and only if the linear system $\dot{x}=\tilde{A}(t)x+\tilde{B}(t)u$ in Problem \eqref{eq:oc6} is controllable over an interval $[0, t_f]$.

% =====================================================
\subsection{The $\b$-coefficients in Theorem \ref{thm:convergence}} 
\label{appd:betas}
The time-varying coefficients $\b_i(\s)$, $i=1,...,9$, in Theorem \ref{thm:convergence} are finite and described below, in which $\s \in [0,t]$ for $0\leq t\leq t_f$:
\begin{align*}
	\beta_1 & = \|(S^{(k)})^T(t)\| \|\Big[(S^{(k)})^T(\s)\Big]^{-1}\|\Big( \|K^{(k)}(\s)\| + \|K^{(k+1)}(\s)\| \Big) \|\Big[S^{(k+1)}(\s)\Big]^{-1} \| \|S^{(k+1)}(t)\|, \\
	\beta_2 & = \|(S^{(k)})^T(t)\| \|\Big[(S^{(k)})^T(\s)\Big]^{-1}\| \|K^{(k)}(\s)\| \|K^{(k+1)}(\s)\| \|\Big[S^{(k+1)}(\s)\Big]^{-1} \| \|S^{(k+1)}(t)\|, \\
	\beta_3 & =  \|\Big[ (S^{(k+1)})^T(t)\Big] ^{-1}\| \|(S^{(k+1)})^T(\s)\| \|x^{(k)}(\sigma)\|, \\
	\beta_4 & =  \|\Big[ (S^{(k+1)})^T(t)\Big] ^{-1}\| \|(S^{(k+1)})^T(\s)\| \| \tilde{B}^{(k-1)}R^{-1}(\tilde{B}^{(k-1)})^T (\sigma) \| \|x^{(k)}(\sigma)\|, \\
	\beta_5 & =  \|\Big[ (S^{(k+1)})^T(t)\Big] ^{-1}\| \|(S^{(k+1)})^T(\s)\| \Big[ \|K^{(k+1)}(\sigma)\| \|x^{(k)}(\sigma)\| + \|S^{(k+1)}\nu^{(k+1)}\| \Big],
\end{align*}
\begin{align*}
	\beta_6 & =  \|\Big[ (S^{(k+1)})^T(t)\Big] ^{-1}\| \|(S^{(k+1)})^T(\s)\|  \| \tilde{B}^{(k-1)}R^{-1}(\tilde{B}^{(k-1)})^T (\sigma) \|, \\
	\beta_7 & = \Big\{ \|(P^{(k+1)})^{-1}\|\Big[ (\|S^{(k+1)}\| + \|S^{(k)}\|)\|\tilde{B}^{(k-1)}R^{-1}(\tilde{B}^{(k-1)})^T(\s)\|\|\nu^{(k)}\|+\|x_0\| \Big] + \|\nu^{(k)}\| \Big\} \\
	&\ \cdot \|\Big[ (S^{(k+1)})^T(t)\Big] ^{-1}\| \|(S^{(k+1)})^T(\s)\|\|(S^{(k)})^T(\s)\|, \\
	\beta_8 & =  \Big\{ \|(P^{(k+1)})^{-1}\|\Big[ (\|S^{(k+1)}\| + \|S^{(k)}\|)\|\tilde{B}^{(k-1)}R^{-1}(\tilde{B}^{(k-1)})^T(\s)\|\|\nu^{(k)}\|+\|x_0\| \Big] + \|\nu^{(k)}\| \Big\} \\
	&\ \cdot \|\Big[ (S^{(k+1)})^T(t)\Big] ^{-1}\| \|(S^{(k+1)})^T(\s)\| \|K^{(k+1)}(\s)\| \|(S^{(k)})^T(\s)\| \\
	&\ + \|(P^{(k+1)})^{-1}\|\|S^{(k+1)}(\s)\| \|S^{(k)}(\s)\|\|\nu^{(k)}\|, \\
	\beta_9 & =  \Big\{ \|(P^{(k+1)})^{-1}\|\Big[ (\|S^{(k+1)}\| + \|S^{(k)}\|)\|\tilde{B}^{(k-1)}R^{-1}(\tilde{B}^{(k-1)})^T(\s)\|\|\nu^{(k)}\|+\|x_0\| \Big] + \|\nu^{(k)}\| \Big\} \\
	&\ \cdot \|\Big[ (S^{(k+1)})^T(t)\Big] ^{-1}\| \|(S^{(k+1)})^T(\s)\|\|\tilde{B}^{(k-1)}R^{-1}(\tilde{B}^{(k-1)})^T(\s)\| \|(S^{(k)})^T(\s)\|.
\end{align*}
\subsection{The entries of the matrix $M$ in \eqref{eq:contraction}} %for the contraction mapping}
\label{appd:M}
Let $p=\sqrt{\sum_{i=1}^n \|G_i \| ^2}$ and $q=\sqrt{\sum_{i,j=1}^n \|H_{ij} \| ^2}$, where $G_i$ and $H_i$ are defined as in \eqref{eq:dA} and \eqref{eq:dB}. The entries of the matrix $M$ in \eqref{eq:contraction} satisfy the relations
\begin{align*}
	m_{11} & \varpropto \beta_3 \Big[ (p+\|x^{(k-1)} \|q)\|K^{(k)} \| + q(\|K^{(k)}\|\|x^{(k)}\|+\|S^{(k)}\nu^{(k)})\|) \Big] + \beta_5 q \Big[\|x^{(k)}\|+\|x^{(k-1)}\| \Big], \\
	m_{12} & \varpropto \beta_3 (p+\|x^{(k-1)} \|q)\|x^{(k-1)} \| + \beta_4, \\
	m_{13} & \varpropto \beta_3 (p+\|x^{(k-1)} \|q) + \beta_6, \\
	m_{21} & \varpropto \beta_1 \Big[ (p+\|x^{(k-1)} \|q)\|K^{(k)} \| + q(\|K^{(k)}\|\|x^{(k)}\|+\|S^{(k)}\nu^{(k)})\|) \Big] + \beta_2 q \Big[\|x^{(k)}\|+\|x^{(k-1)}\| \Big], \\
	m_{22} & \varpropto \beta_1 (p+\|x^{(k-1)} \|q)\|x^{(k-1)} \|, \\
	m_{23} & \varpropto \beta_1 (p+\|x^{(k-1)} \|q), \\
	m_{31} & \varpropto \beta_7 \Big[ (p+\|x^{(k-1)} \|q)\|K^{(k)} \| + q(\|K^{(k)}\|\|x^{(k)}\|+\|S^{(k)}\nu^{(k)})\|) \Big] + \beta_8 q \Big[\|x^{(k)}\|+\|x^{(k-1)}\| \Big], \\
	m_{32} & \varpropto \beta_7 (p+\|x^{(k-1)} \|q)\|x^{(k-1)} \| + \beta_9, \\
	m_{33} & \varpropto \beta_7 (p+\|x^{(k-1)} \|q),
\end{align*}
where $\varpropto$ denotes proportionality.
Because both $p$ and $q$ are related to $R^{-1}$, they can be made sufficiently small by choosing large enough $R$, for example, $R=\g I$ with $\g\gg 1$, where $I\in\mathbb{R}^{m\times m}$ is the identify matrix. In addition, for $m_{12}$, $m_{13}$, and $m_{32}$, the coefficients $\beta_4$, $\beta_6$, and $\beta_9$ defined in Appendix \ref{appd:betas} involve the factor $\|\tilde{B}^{(k-1)}R^{-1}(\tilde{B}^{(k-1)})^T(\s)\|$, $\s\in [0,t_f]$, which can also be made arbitrary small by adjusting $R$, more specifically, by choosing $R$ with large eigenvalues. Thus, each entry of $M$ can be made sufficiently small through the choice of $R$.

% =================================================
\subsection{Singular value expansion for compact operators}
\label{appd:linear_ensemble}

% ============= Theorem 8.3 (Singular Value Expansion) ==============
\begin{theorem}[Singular value expansion \cite{Gohberg}]
	\label{thm:singular}
	Let $Y$ and $Z$ be Hilbert spaces, $K:Y\rightarrow Z$ be a compact operator and $\{(\sigma_n,\mu_n,\nu_n)\ |\ n\in\Delta\}$ be a singular system for $K$. Then
	$$Ky=\sum_{n\in \Delta }\sigma_n \langle  y,\mu_n \rangle \nu_n,\quad K^* z=\sum_{n\in \Delta } \sigma_n \langle z,\nu_n \rangle \mu_n,$$
	for all $y\in Y$, $z\in Z$. In particular, if $K_n y=\sum_{j=1}^{n}\sigma_j \langle y,\mu_j \rangle \nu_j$ for $y\in Y$, and $K$ is of infinite rank, namely, $\Delta=\mathbb{N}$, then $\|K-K_n\|\leq\sup_{j>n}\sigma_j\rightarrow 0$ as $n\rightarrow\infty$.
\end{theorem}

%%%%%%%%%%%%%%%%%%%%%%%%%%%%%%%%%%%%%%%%%%%%%%%%%%%%%%%%%%%%%%%%%%%%%%%%%%
\section{Bibliography}
\bibliographystyle{ieeetr}
\bibliography{Iterative}

\end{document}